\documentclass[12 pt]{article}
\usepackage{amssymb,amsmath,amsfonts,amsthm}

\theoremstyle{plain}
\newtheorem{thm}{Theorem}
\newtheorem{lma}{Lemma}
\newtheorem{rmk}{Remark}
\newtheorem{prop}{Proposition}
\newtheorem{cor}{Corollary}

\theoremstyle{definition}
\newtheorem{defn}{Definition}

\newtheorem{nota}{Notation}

\begin{document}

\begin{center}
\textbf{\Large A Bernstein type inequality and moderate deviations for
weakly dependent sequences}\vskip15pt

Florence Merlev\`{e}de $^{a}$, Magda Peligrad $^{b}$ \footnote{%
Supported in part by a Charles Phelps Taft Memorial Fund grant, and NSA\
grants, H98230-07-1-0016 and H98230-09-1-0005} \textit{and\/} Emmanuel Rio $%
^{c}$ \footnote{%
Supported in part by Centre INRIA Bordeaux Sud-Ouest $\&$ Institut de
Math\'ematiques de Bordeaux}
\end{center}

\vskip10pt $^a$ Universit\'e Paris Est, Laboratoire de math\'ematiques, UMR
8050 CNRS, B\^atiment Copernic, 5 Boulevard Descartes, 77435
Champs-Sur-Marne, FRANCE. E-mail: florence.merlevede@univ-mlv.fr\newline
\newline
$^b$ Department of Mathematical Sciences, University of Cincinnati, PO Box
210025, Cincinnati, Oh 45221-0025, USA. Email: magda.peligrad@uc.edu \newline
\newline
$^c$ Universit\'e de Versailles, Laboratoire de math\'ematiques, UMR 8100
CNRS, B\^atiment Fermat, 45 Avenue des Etats-Unis, 78035 Versailles, FRANCE.
E-mail: rio@math.uvsq.fr \vskip10pt



\begin{center}
\textbf{Abstract}\vskip10pt
\end{center}

In this paper we present a tail inequality for the maximum of partial sums
of a weakly dependent sequence of random variables that is not necessarily
bounded. The class considered includes geometrically and subgeometrically
strongly mixing sequences. The result is then used to derive asymptotic
moderate deviation results. Applications include classes of Markov chains,
functions of linear processes with absolutely regular innovations and ARCH
models.

\section{Introduction}

\setcounter{equation}{0}

Let us consider a sequence $X_{1},X_{2},\ldots $ of real valued random
variables. The aim of this paper is to present nonasymptotic tail
inequalities for $S_{n}=X_{1}+X_{2}+\cdots +X_{n}$ and to use them to derive
moderate deviations principles.

\smallskip For independent and centered random variables $X_{1},X_{2},\ldots 
$, one of the main tools to get an upper bound for the large and moderate
deviations principles is the so-called Bernstein inequalities. We first
recall the Bernstein inequality for random variables satisfying Condition (%
\ref{b1lt}) below. Suppose that the random variables $X_{1},X_{2},\ldots $
satisfy 
\begin{equation}
\log \mathbb{E}\exp (tX_{i})\leq \frac{\sigma _{i}^{2}t^{2}}{2(1-tM)}\ 
\hbox{ for
positive constants }\ \sigma _{i}\ \hbox{ and }M,  \label{b1lt}
\end{equation}%
for any $t$ in $[0,1/M[$. Set $V_{n}=\sigma _{1}^{2}+\sigma _{2}^{2}+\cdots
+\sigma _{n}^{2}$. Then 
\begin{equation*}
\mathbb{P}(S_{n}\geq \sqrt{2V_{n}x}+Mx)\leq \exp (-x).
\end{equation*}%
When the random variables $X_{1},X_{2},\ldots $ are uniformly bounded by $M$
then (\ref{b1lt}) holds with $\sigma _{i}^{2}=\mathrm{Var}X_{i}$, and the
above inequality implies the usual Bernstein inequality 
\begin{equation}
\mathbb{P}(S_{n}\geq y)\leq \exp \Bigl(-y^{2}(2V_{n}+2yM)^{-1}\Bigr).
\label{B1}
\end{equation}%
Assume now that the random variables $X_{1},X_{2},\ldots $ satisfy the
following weaker tail condition: for some $\gamma $ in $]0,1[$ and any positive $t$, 
$\sup_i \mathbb{P}(X_{i}\geq t)\leq \exp (1-(t/M)^{\gamma })$. Then, by the proof of Corollary 5.1 in
Borovkov (2000-a) we infer that 
\begin{equation}
\mathbb{P}(|S_{n}|\geq y)\leq 2\exp \Bigl(-c_{1}y^{2}/V_{n}\Bigr)+n\exp %
\Bigl(-c_{2}(y/M)^{\gamma }\Bigr)\,,  \label{B2}
\end{equation}%
where $c_{1}$ and $c_{2}$ are positive constants ($c_{2}$ depends on $\gamma 
$). More precise results for large and moderate deviations of sums of
independent random variables with \textsl{semiexponential\/} tails may be
found in Borovkov (2000-b).

In our terminology the moderate deviations principle (MDP) stays for the
following type of asymptotic behavior:

\smallskip

\begin{defn}
We say that the MDP holds for a sequence $(T_{n})_{n}$ of random variables
with the speed $a_{n}\rightarrow 0$ and rate function $I(t)$ if for each $A$
Borelian, 
\begin{eqnarray}
-\inf_{t\in A^{o}}I(t) &\leq &\liminf_{n}a_{n}\log \mathbb{P}(\sqrt{a_{n}}%
T_{n}\in A)  \notag \\
&\leq &\limsup_{n}a_{n}\log \mathbb{P}(\sqrt{a_{n}}T_{n}\in A)\leq
-\inf_{t\in \bar{A}}I(t)\,,  \label{mdpdef}
\end{eqnarray}%
where $\bar{A}$ denotes the closure of $A$ and $A^{o}$ the interior of $A$.
\end{defn}

Our interest is to extend the above inequalities to strongly mixing
sequences of random variables and to study the MDP for $%
(S_{n}/stdev(S_{n}))_{n}$. In order to cover a larger class of examples we
shall also consider less restrictive coefficients of weak dependence, such
as the $\tau $-mixing coefficients defined in Dedecker and Prieur (2004)
(see Section \ref{sectionMR} for the definition of these coefficients).

Let $X_{1},X_{2},\ldots $ be a strongly mixing sequence of real-valued and
centered random variables. Assume that there exist a positive constant $%
\gamma _{1}$ and a positive $c$ such that the strong mixing coefficients of
the sequence satisfy 
\begin{equation}
\alpha (n)\leq \exp (-cn^{\gamma _{1}})\ \hbox{ for any positive integer }%
n\,,  \label{H1}
\end{equation}%
and there is a constant $\gamma _{2}$ in $]0,+\infty ]$ such that 
\begin{equation}
\sup_{i>0}\mathbb{P}(|X_{i}|>t)\leq \exp (1-t^{\gamma _{2}})\ 
\hbox{ for any
positive }t  \label{H2}
\end{equation}%
(when $\gamma _{2}=+\infty $ (\ref{H2}) means that $\Vert X_{i}\Vert
_{\infty }\leq 1$ for any positive $i$).

Obtaining exponential bounds for this case is a challenging problem. One of
the available tools in the literature is Theorem 6.2 in Rio (2000), which is
a Fuk-Nagaev type inequality, that provides the inequality below. Let $%
\gamma $ be defined by $1/\gamma =(1/\gamma _{1})+(1/\gamma _{2})$. For any
positive $\lambda $ and any $r\geq 1$, 
\begin{equation}
\mathbb{P}(\sup_{k\in \lbrack 1,n]}|S_{k}|\geq 4\lambda )\leq 4\Bigl(1+\frac{%
\lambda ^{2}}{rnV}\Bigr)^{-r/2}+4Cn\lambda ^{-1}\exp \Bigl(-c(\lambda
/r)^{\gamma }\Bigr),  \label{RFN}
\end{equation}%
where 
\begin{equation*}
V=\sup_{i>0}\Bigl(\mathbb{E}(X_{i}^{2})+2\sum_{j>i}|\mathbb{E}(X_{i}X_{j})|%
\Bigr).
\end{equation*}%
Selecting in (\ref{RFN}) $r=\lambda ^{2}/(nV)$ leads to 
\begin{equation*}
\mathbb{P}(\sup_{k\in \lbrack 1,n]}|S_{k}|\geq 4\lambda )\leq 4\exp \Bigl(-%
\frac{\lambda ^{2}\log 2}{2nV}\Bigr)+4Cn\lambda ^{-1}\exp \Bigl(%
-c(nV/\lambda )^{\gamma }\Bigr)
\end{equation*}%
for any $\lambda \geq (nV)^{1/2}$. The above inequality gives a subgaussian
bound, provided that 
\begin{equation*}
(nV/\lambda )^{\gamma }\geq \lambda ^{2}/(nV)+\log (n/\lambda ),
\end{equation*}%
which holds if $\lambda \ll (nV)^{(\gamma +1)/(\gamma +2)}$ (here and below $%
\ll $ replaces the 
symbol $o$). Hence (\ref{RFN}) is useful to study the probability of
moderate deviation $\mathbb{P}(|S_{n}|\geq t\sqrt{n/a_{n}})$ provided $%
a_{n}\gg n^{-\gamma /(\gamma +2)}$. For $\gamma =1$ this leads to $a_{n}\gg
n^{-1/3}$. For bounded random variables and geometric mixing rates (in that
case $\gamma =1$), Proposition 13 in Merlev\`{e}de and Peligrad (2009)
provides the MDP under the improved condition $a_{n}\gg n^{-1/2}$. We will
prove in this paper that this condition is still suboptimal from the point
of view of moderate deviation.

\smallskip For stationary geometrically mixing (absolutely regular) Markov
chains, and bounded functions $f$ (here $\gamma =1$), Theorem 6 in Adamczak
(2008) provides a Bernstein's type inequality for $%
S_{n}(f)=f(X_{1})+f(X_{2})+\cdots +f(X_{n})$. Under the centering condition $%
\mathbb{E}(f(X_{1}))=0$, he proves that 
\begin{equation}
\mathbb{P}(|S_{n}(f)|\geq \lambda )\leq C\exp \Bigl(-\frac{1}{C}\min \Bigl(%
\frac{\lambda ^{2}}{n\sigma ^{2}},\frac{\lambda }{\log n}\Bigr)\Bigr),
\label{Adam}
\end{equation}%
where $\sigma ^{2}=\lim_{n}n^{-1}\mathrm{Var}S_{n}(f)$ (here we take $m=1$
in his condition (14) on the small set). Inequality (\ref{Adam}) provides
exponential tightness for $S_{n}(f)/\sqrt{n}$ with rate $a_{n}$ as soon as $%
a_{n}\gg n^{-1}(\log n)^{2}$, which is weaker than the above conditions.
Still in the context of Markov chains, we point out the recent Fuk-Nagaev
type inequality obtained by Bertail and Cl\'emen\c{c}on (2008). However for
stationary subgeometrically mixing Markov chains, their inequality does not
lead to the optimal rate which can be expected in view of the results
obtained by Djellout and Guillin (2001).

\smallskip To our knowledge, Inequality (\ref{Adam}) has not been extended
yet to the case $\gamma <1$, even for the case of bounded functions $f$ and
absolutely regular Markov chains. In this paper we improve inequality (\ref%
{RFN}) in the case $\gamma <1 $ and then derive moderate deviations
principles from this new inequality under the minimal condition $%
a_{n}n^{\gamma /(2-\gamma )}\rightarrow \infty $. The main tool is an
extension of inequality (\ref{B2}) to dependent sequences. We shall prove
that, for $\alpha $-mixing or $\tau $-mixing sequences satisfying (\ref{H1})
and (\ref{H2}) for $\gamma <1$, there exists a positive $\eta $ such that,
for $n\geq 4$ and $\lambda \geq C(\log n)^{\eta }$ 
\begin{equation}
\mathbb{P}( \sup_{j\leq n}|S_{j}| \geq \lambda )\leq (n+1) \exp
(-\lambda ^{\gamma } / C_1) +  \exp (-\lambda ^{2}/ (C_2 + C_2nV) ) ,  \label{Boro}
\end{equation}%
where $C$, $C_1$ and $C_2$ are positive constants depending on $c$, $\gamma
_{1}$ and $\gamma _{2}$ and $V$ is some constant  (which differs from the
constant $V$ in (\ref{RFN}) in the unbounded case), depending on the
covariance properties of truncated random variables built from the initial
sequence. In order to define precisely $V$ we need to introduce truncation
functions $\varphi _{M}$.

\begin{nota}
For any positive $M$ let the function $\varphi_M$ be defined by $\varphi_M
(x) = (x \wedge M) \vee(-M) $.
\end{nota}

With this notation, (\ref{Boro}) holds with 
\begin{equation}
V=\sup_{M\geq 1}\sup_{i>0}\Bigl (\mathrm{Var} (\varphi
_{M}(X_{i}))+2\sum_{j>i}|\mathrm{Cov}(\varphi _{M}(X_{i}),\varphi
_{M}(X_{j}))|\Bigr).  \label{VV}
\end{equation}%
To prove (\ref{Boro}) we use a variety of techniques and new ideas, ranging
from the big and small blocks argument based on a Cantor-type construction,
diadic induction, adaptive truncation along with coupling arguments. In a
forthcoming paper, we will study the case $\gamma _{1}=1$ and $\gamma
_{2}=\infty $. We now give more definitions and precise results.

\section{Main results}

\label{sectionMR}

\setcounter{equation}{0}

\quad We first define the dependence coefficients that we consider in this
paper.

\quad For any real random variable $X$ in ${\mathbb{L}}^{1}$ and any $\sigma 
$-algebra $\mathcal{M}$ of $\mathcal{A}$, let ${\mathbb{P}}_{X|\mathcal{M}}$
be a conditional distribution of $X$ given ${\mathcal{M}}$ and let ${\mathbb{%
\ P}}_{X}$ be the distribution of $X$. We consider the coefficient $\tau ( 
\mathcal{M},X)$ of weak dependence (Dedecker and Prieur, 2004) which is
defined by 
\begin{equation}
\tau (\mathcal{M},X)=\Big \|\sup_{f\in \Lambda _{1}(\mathbb{R})}\Bigr|\int
f(x)\mathbb{P}_{X|\mathcal{M}}(dx)-\int f(x)\mathbb{P}_{X}(dx)\Big |\Big \|%
_{1}\, ,  \label{deftau1}
\end{equation}%
where $\Lambda _{1}(\mathbb{R})$ is the set of $1$-Lipschitz functions from $%
\mathbb{R}$ to $\mathbb{R}$.

The coefficient $\tau $ has the following coupling property: If $\Omega $ is
rich enough then the coefficient $\tau (\mathcal{M},X)$ is the infimum of $%
\Vert X-X^{\ast }\Vert _{1}$ where $X^{\ast }$ is independent of $\mathcal{M}
$ and distributed as $X$ (see Lemma 5 in Dedecker and Prieur (2004)). This
coupling property allows to relate the coefficient $\tau $ to the strong
mixing coefficient Rosenblatt (1956) defined by 
\begin{equation*}
\alpha (\mathcal{M},\sigma (X))=\sup_{A\in \mathcal{M},B\in \sigma (X)}|{%
\mathbb{P}}(A\cap B)-{\mathbb{P}}(A){\mathbb{P}}(B)|\,,
\end{equation*}%
as shown in Rio (2000, p. 161) for the bounded case, and by Peligrad (2002)
for the unbounded case. For equivalent definitions of the strong mixing
coefficient we refer for instance to Bradley (2007, Lemma 4.3 and Theorem
4.4).

\medskip

\quad If $Y$ is a random variable with values in $\mathbb{R}^{k}$, the
coupling coefficient $\tau$ is defined as follows: If $Y\in {\mathbb{L}}^{1}(%
\mathbb{R}^{k})$, 
\begin{equation}
\tau (\mathcal{M},Y)=\sup \{\tau (\mathcal{M},f(Y)),f\in \Lambda _{1}(%
\mathbb{R}^{k})\}\,,  \label{deftau1}
\end{equation}%
where $\Lambda _{1}(\mathbb{R}^{k})$ is the set of $1$-Lipschitz functions
from $\mathbb{R}^{k}$ to $\mathbb{R}$.

\medskip The $\tau $-mixing coefficients $\tau _{X}(i)=\tau (i)$ of a
sequence $(X_{i})_{i\in \mathbb{Z}}$ of real-valued random variables are
defined by 
\begin{equation}
\tau _{k}(i)=\max_{1\leq \ell \leq k}\frac{1}{\ell }\sup \Big \{\tau (%
\mathcal{M}_{p},(X_{j_{1}},\cdots ,X_{j_{\ell }})),\,p+i\leq j_{1}<\cdots
<j_{\ell }\Big \}\text{ and }\tau (i)=\sup_{k\geq 0}\tau _{k}(i)\,,
\label{deftau2}
\end{equation}%
where $\mathcal{M}_{p}=\sigma (X_{j},j\leq p)$ and the above supremum is
taken over $p$ and $(j_{1},\ldots j_{\ell})$. Recall that the strong mixing
coefficients $\alpha (i)$ are defined by: 
\begin{equation*}
\ {\alpha }(i)=\sup_{p\in \mathbb{Z}}\alpha (\mathcal{M}_{p},\sigma
(X_{j},j\geq i+p))\,.
\end{equation*}%
Define now the function $Q_{|Y|}$ by $Q_{|Y|}(u)=\inf \{t>0,\mathbb{P}%
(|Y|>t)\leq u\}$ for $u$ in $]0,1]$. To compare the $\tau $-mixing
coefficient with the strong mixing coefficient, let us mention that, by
Lemma 7 in Dedecker and Prieur (2004), 
\begin{equation}
\tau (i)\leq 2\int_{0}^{2\alpha (i)}Q(u)du\,,\hbox{ where }Q=\sup_{k\in {%
\mathbb{Z}}}Q_{|X_{k}|}.  \label{comptaualpha}
\end{equation}

\smallskip Let $(X_{j})_{j\in \mathbb{Z}}$ be a sequence of centered real
valued random variables and let $\tau (i)$ be defined by (\ref{deftau2}).
Let $\tau (x)=\tau ([x])$ (square brackets denoting the integer part).
Throughout, we assume that there exist positive constants $\gamma _{1}$ and $%
\gamma _{2}$ such that 
\begin{equation}
\tau (x)\leq \exp (-cx^{\gamma _{1}})\ \hbox{ for any }x\geq 1\,,
\label{hypoalpha}
\end{equation}%
where $c>0$ and for any positive $t$, 
\begin{equation}
\sup_{k>0}\mathbb{P}(|X_{k}|>t)\leq \exp (1-t^{\gamma _{2}}):=H(t)\,.
\label{hypoQ}
\end{equation}%
Suppose furthermore that 
\begin{equation}
\gamma <1\text{ where $\gamma $ is defined by $1/\gamma =1/\gamma
_{1}+1/\gamma _{2}$}\,.  \label{hypogamma}
\end{equation}

\begin{thm}
\label{BTinegacont} Let $(X_{j})_{j\in \mathbb{Z}}$ be a sequence of
centered real valued random variables and let $V$ be defined by (\ref{VV}).
Assume that (\ref{hypoalpha}), (\ref{hypoQ}) and (\ref{hypogamma}) are
satisfied. Then $V$ is finite and, for any $n\geq 4$, there exist positive
constants $C_{1}$, $C_{2}$, $C_{3}$ and $C_{4}$ depending only on $c $, $%
\gamma $ and $\gamma _{1}$ such that, for any positive $x$, 
\begin{equation*}
\mathbb{P}\Bigl(\sup_{j\leq n}|S_{j}|\geq x\Bigr)\leq n\exp \Bigl(-\frac{%
x^{\gamma }}{C_{1}}\Bigr )+\exp \Bigl(-\frac{x^{2}}{C_{2}(1+nV)}\Bigr)+\exp %
\Bigl(-\frac{x^{2}}{C_{3}n}\exp \Bigr(\frac{x^{\gamma (1-\gamma )}}{%
C_{4}(\log x)^{\gamma }}\Bigr)\Bigr)\,.
\end{equation*}
\end{thm}

\begin{rmk}
Let us mention that if the sequence $(X_{j})_{j\in \mathbb{Z}}$ satisfies (%
\ref{hypoQ}) and is strongly mixing with strong mixing coefficients
satisfying (\ref{H1}), then, from (\ref{comptaualpha}), (\ref{hypoalpha}) is
satisfied (with an other constant), and Theorem \ref{BTinegacont} applies.
\end{rmk}

\begin{rmk}
\label{remexpmom} If $\mathbb{E} \exp ( |X_i|^{\gamma_2})) \leq K$ for any
positive $i$, then setting $C = 1 \vee \log K $, we notice that the process $%
(C^{-1/\gamma_2}X_i)_{i \in {\mathbb{Z}}}$ satisfies (\ref{hypoQ}).
\end{rmk}

\begin{rmk}
\label{remv2} If $(X_{i})_{i\in \mathbb{Z} }$ satisfies (\ref{hypoalpha})
and (\ref{hypoQ}), then 
\begin{eqnarray*}
V &\leq &\sup_{i>0} \Bigl( \mathbb{E} (X_i^2) + 4 \sum_{k >0 }
\int_0^{\tau(k)/2} Q_{|X_i|} (G(v)) dv \Bigr) \\
& = & \sup_{i>0} \Bigl( \mathbb{E} (X_i^2) + 4 \sum_{k >0 }
\int_0^{G(\tau(k)/2)} Q_{|X_i|} (u) Q(u) du \Bigr) ,
\end{eqnarray*}
where $G$ is the inverse function of $x \mapsto \int_0^xQ(u) du$ (see
Section \ref{prremv2} for a proof). Here the random variables do not need to
be centered. Note also that, in the strong mixing case, using (\ref%
{comptaualpha}), we have $G(\tau(k)/2) \leq 2 \alpha (k)$.
\end{rmk}

This result is the main tool to derive the MDP below.

\begin{thm}
\label{thmMDPsubgeo11} Let $(X_{i})_{i\in \mathbb{Z}}$ be a sequence of
random variables as in Theorem \ref{BTinegacont} and let $%
S_{n}=\sum_{i=1}^{n}X_{i}$ and $\sigma _{n}^{2}=\mathrm{Var}S_{n}$. Assume
in addition that $\lim \inf_{n\rightarrow \infty }\sigma _{n}^{2}/n>0$. Then
for all positive sequences $a_{n}$ with $a_{n}\rightarrow 0$ and $%
a_{n}n^{\gamma /(2-\gamma )}\rightarrow \infty $, $\{\sigma _{n}^{-1}S_{n}\}$
satisfies (\ref{mdpdef}) with the good rate function $I(t)=t^{2}/2$.
\end{thm}

If we impose a stronger degree of stationarity we obtain the following
corollary.

\begin{cor}
\label{thmMDPsubgeo} Let $(X_{i})_{i\in \mathbb{Z}}$ be a second order
stationary sequence of centered real valued random variables. Assume that (%
\ref{hypoalpha}), (\ref{hypoQ}) and (\ref{hypogamma}) are satisfied. Let $%
S_{n}=\sum_{i=1}^{n}X_{i}$ and $\sigma _{n}^{2}=\mathrm{Var}S_{n}$. Assume
in addition that $\sigma _{n}^{2}\rightarrow \infty $. Then $%
\lim_{n\rightarrow \infty }\sigma _{n}^{2}/n=\sigma ^{2}>0$, and for all
positive sequences $a_{n}$ with $a_{n}\rightarrow 0$ and $a_{n}n^{\gamma
/(2-\gamma )}\rightarrow \infty $, $\{n^{-1/2}S_{n}\}$ satisfies (\ref%
{mdpdef}) with the good rate function $I(t)=t^{2}/(2\sigma ^{2})$.
\end{cor}

\subsection{Applications}

\subsubsection{Instantaneous functions of absolutely regular processes}

Let $(Y_{j})_{j\in \mathbb{Z}}$ be a strictly stationary sequence of random
variables with values in a Polish space $E$, and let $f$ be a measurable
function from $E$ to ${\mathbb{R}}$. Set $X_{j}=f(Y_{j})$. Consider now the
case where the sequence $(Y_{k})_{k\in \mathbb{Z}}$ is absolutely regular
(or $\beta $-mixing) in the sense of Rozanov and Volkonskii (1959). Setting $%
\mathcal{F}_{0}=\sigma (Y_{i},i\leq 0)$ and $\mathcal{G}_{k}=\sigma
(Y_{i},i\geq k)$, this means that 
\begin{equation*}
\beta (k)=\beta (\mathcal{F}_{0},\mathcal{G}_{k})\rightarrow 0\,,\text{ as }%
k\rightarrow \infty \,,
\end{equation*}%
with $\beta ({\mathcal{A}},{\mathcal{B}})=\frac{1}{2}\sup \{\sum_{i\in
I}\sum_{j\in J}|\mathbb{P}(A_{i}\cap B_{j})-\mathbb{P}(A_{i})\mathbb{P}%
(B_{j})|\}$, the maximum being taken over all finite partitions $%
(A_{i})_{i\in I}$ and $(B_{i})_{i\in J}$ of $\Omega $ respectively with
elements in ${\mathcal{A}}$ and ${\mathcal{B}}$. If we assume that 
\begin{equation}
\beta (n)\leq \exp (-cn^{\gamma _{1}})\text{ for any positive }n,
\label{hypobeta}
\end{equation}%
where $c>0$ and $\gamma _{1}>0$, and that the random variables $X_{j}$ are
centered and satisfy (\ref{hypoQ}) for some positive $\gamma _{2}$ such that 
$1/\gamma =1/\gamma _{1}+1/\gamma _{2}>1$, then Theorem \ref{BTinegacont}
and Corollary \ref{thmMDPsubgeo} apply to the sequence $(X_{j})_{j\in 
\mathbb{Z}}$. Furthermore, as shown in Viennet (1997), by Delyon's (1990)
covariance inequality, 
\begin{equation*}
V\leq \mathbb{E}(f^{2}(X_{0}))+4\sum_{k>0}\mathbb{E}(B_{k}f^{2}(X_{0})),
\end{equation*}%
for some sequence $(B_{k})_{k>0}$ of random variables with values in $[0,1]$
satisfying $\mathbb{E}(B_{k})\leq \beta (k)$ (see Rio (2000, Section 1.6)
for more details).

We now give an example where $(Y_{j})_{j\in \mathbb{Z}}$ satisfies (\ref%
{hypobeta}). Let $(Y_{j})_{j\geq 0}$ be an $E$-valued irreducible ergodic
and stationary Markov chain with a transition probability $P$ having a
unique invariant probability measure $\pi $ (by Kolmogorov extension Theorem
one can complete $(Y_{j})_{j\geq 0}$ to a sequence $(Y_{j})_{j\in \mathbb{Z}%
} $). Assume furthermore that the chain has an atom, that is there exists $%
A\subset E$ with $\pi (A)>0$ and $\nu $ a probability measure such that $%
P(x,\cdot )=\nu (\cdot )$ for any $x$ in $A$. If 
\begin{equation}
\text{ there exists }\delta >0\text{ and }\gamma _{1}>0\text{ such that }%
\mathbb{E}_{\nu }(\exp (\delta \tau ^{\gamma _{1}}))<\infty \,,
\label{condsg}
\end{equation}%
where $\tau =\inf \{n\geq 0;\,Y_{n}\in A\}$, then the $\beta $-mixing
coefficients of the sequence $(Y_{j})_{j\geq 0}$ satisfy (\ref{hypobeta})
with the same $\gamma _{1}$ (see Proposition 9.6 and Corollary 9.1 in Rio
(2000) for more details). Suppose that $\pi (f)=0$. Then the results apply
to $(X_{j})_{j\geq 0}$ as soon as $f$ satisfies 
\begin{equation*}
\pi (|f|>t)\leq \exp (1-t^{\gamma _{2}})\ \text{ for any positive }t\,.
\end{equation*}%
Compared to the results obtained by de Acosta (1997) and Chen and de Acosta
(1998) for geometrically ergodic Markov chains, and by Djellout and Guillin
(2001) for subgeometrically ergodic Markov chains, we do not require here
the function $f$ to be bounded.

\subsubsection{Functions of linear processes with absolutely regular
innovations}

\label{Sexample1}

Let $f$ be a 1-Lipshitz function. We consider here the case where 
\begin{equation*}
X_{n}=f\bigl(\sum_{j\geq 0}a_{j}\xi _{n-j}\bigr)-\mathbb{E}f\bigl(%
\sum_{j\geq 0}a_{j}\xi _{n-j}\bigr)\,,
\end{equation*}%
where $A=\sum_{j\geq 0}|a_{j}|<\infty $ and $(\xi _{i})_{i\in {\mathbb{Z}}}$
is a strictly stationary sequence of real-valued random variables which is
absolutely regular in the sense of Rozanov and Volkonskii.

Let $\mathcal{F}_{0}=\sigma (\xi _{i},i\leq 0)$ and $\mathcal{G}_{k}=\sigma (%
{\xi }_{i},i\geq k)$. According to Section 3.1 in Dedecker and Merlev\`{e}de
(2006), if the innovations $(\xi _{i})_{i\in {\mathbb{Z}}}$ are in ${\mathbb{%
L}}^{2}$, the following bound holds for the $\tau $-mixing coefficient
associated to the sequence $(X_{i})_{i\in {\mathbb{Z}}}$: 
\begin{equation*}
\tau (i)\leq 2\Vert \xi _{0}\Vert _{1}\sum_{j\geq i}|a_{j}|+4\Vert \xi
_{0}\Vert _{2}\sum_{j=0}^{i-1}|a_{j}|\,\beta _{\xi }^{1/2}(i-j)\,.
\end{equation*}%
Assume that there exists $\gamma _{1}>0$ and $c^{\prime }>0$ such that, for
any positive integer $k$, 
\begin{equation*}
a_{k}\leq \exp (-c^{\prime }k^{\gamma _{1}})\text{ and }\beta _{\xi }(k)\leq
\exp (-c^{\prime }k^{\gamma _{1}})\,.
\end{equation*}%
Then the $\tau $-mixing coefficients of $(X_{j})_{j\in \mathbb{Z}}$ satisfy (%
\ref{hypoalpha}). Let us now focus on the tails of the random variables $%
X_{i}$. Assume that $(\xi _{i})_{i\in {\mathbb{Z}}}$ satisfies (\ref{hypoQ}%
). Define the convex functions $\psi _{\eta }$ for $\eta >0$ in the
following way: $\psi _{\eta }(-x)=\psi _{\eta }(x)$, and for any $x\geq 0$, 
\begin{equation*}
\psi _{\eta }(x)=\exp (x^{\eta })-1\hbox{ for }\eta \geq 1\ \hbox{ and }\
\psi _{\eta }(x)=\int_{0}^{x}\exp (u^{\eta })du\hbox{ for }\eta \in ]0,1].
\end{equation*}%
Let $\Vert \,.\,\Vert _{\psi _{\eta }}$ be the usual corresponding Orlicz
norm. Since the function $f$ is 1-Lipshitz, we get that $\Vert X_{0}\Vert
_{\psi _{\gamma _{2}}}\leq 2A\Vert \xi _{0}\Vert _{\psi _{\gamma _{2}}}$.
Next, if $(\xi _{i})_{i\in {\mathbb{Z}}}$ satisfies (\ref{hypoQ}), then $%
\Vert \xi _{0}\Vert _{\psi _{\gamma _{2}}}<\infty $. Furthermore, it can
easily be proven that, if $\Vert Y\Vert _{\psi _{\eta }}\leq 1$, then $%
\mathbb{P}(|Y|>t)\leq \exp (1-t^{\eta })$ for any positive $t$. Hence,
setting $C=2A\Vert \xi _{0}\Vert _{\psi _{\gamma _{2}}}$, we get that $%
(X_{i}/C)_{i\in {\mathbb{Z}}}$ satisfies (\ref{hypoQ}) with the same
parameter $\gamma _{2}$, and therefore the conclusions of Theorem \ref%
{BTinegacont} and Corollary \ref{thmMDPsubgeo} hold with $\gamma $ defined
by $1/\gamma =1/\gamma _{1}+1/\gamma _{2}$, provided that $\gamma <1$.

This example shows that our results hold for processes that are not
necessarily strongly mixing. Recall that, in the case where $a_{i}=2^{-i-1}$
and the innovations are iid with law ${\mathcal{B}}(1/2)$, the process fails
to be strongly mixing in the sense of Rosenblatt.

\subsubsection{ARCH($\infty$) models}

\label{Sexample2} Let $(\eta _{t})_{t\in \mathbb{Z}}$ be an iid sequence of
zero mean real random variables such that $\Vert \eta _{0}\Vert _{\infty
}\leq 1$. We consider the following ARCH($\infty $) model described by
Giraitis \textit{et al.} (2000): 
\begin{equation}
Y_{t}=\sigma _{t}\eta _{t}\,,\hbox{ where }\sigma _{t}^{2}=a+\sum_{j\geq
1}a_{j}Y_{t-j}^{2}\,,  \label{defarch}
\end{equation}%
where $a\geq 0$, $a_{j}\geq 0$ and $\sum_{j\geq 1}a_{j}<1$. Such models are
encountered, when the volatility $(\sigma _{t}^{2})_{t\in {\mathbb{Z}}}$ is
unobserved. In that case, the process of interest is $(Y_{t}^{2})_{t\in {%
\mathbb{Z}}}$. Under the above conditions, there exists a unique stationary
solution that satisfies 
\begin{equation*}
\Vert Y_{0}\Vert _{\infty }^{2}\leq a+a\sum_{\ell \geq 1}\big (\sum_{j\geq
1}a_{j}\big )^{\ell }=M<\infty \,.
\end{equation*}

Set now $X_{j}=(2M)^{-1}(Y_{j}^{2}-\mathbb{E}(Y_{j}^{2}))$. Then the
sequence $(X_{j})_{j\in \mathbb{Z}}$ satisfies (\ref{hypoQ}) with $\gamma
_{2}=\infty $. If we assume in addition that $a_{j}=O(b^{j})$ for some $b<1$%
, then, according to Proposition 5.1 (and its proof) in Comte \textit{et al.}
(2008), the $\tau $-mixing coefficients of $(X_{j})_{j\in \mathbb{Z}}$
satisfy (\ref{hypoalpha}) with $\gamma _{1}=1/2$. Hence in this case, the
sequence $(X_{j})_{j\in \mathbb{Z}}$ satisfies both the conclusions of
Theorem \ref{BTinegacont} and of Corollary \ref{thmMDPsubgeo} with $\gamma
=1/2$.

\section{Proofs}

\setcounter{equation}{0}

\subsection{Some auxiliary results}

The aim of this section is essentially to give suitable bounds for the
Laplace transform of 
\begin{equation}  \label{Bp1}
S (K) = \sum_{i \in K} X_i \, ,
\end{equation}
where $K$ is a finite set of integers.

\medskip

\begin{equation}
c_{0}=(2(2^{1/\gamma }-1))^{-1}(2^{(1-\gamma )/\gamma }-1)\,,\ c_{1}=\min
(c^{1/\gamma _{1}}c_{0}/4,2^{-1/\gamma })\,,  \label{defc0}
\end{equation}%
\begin{equation}
c_{2}=2^{-(1+2\gamma _{1}/\gamma )}c_{1}^{\gamma
_{1}}\,\,,\,\,c_{3}=2^{-\gamma _{1}/\gamma }\,,\text{ and }\kappa =\min \big
(c_{2},c_{3}\big )\,.  \label{defkappa}
\end{equation}

\begin{prop}
\label{propinter2} Let $(X_{j})_{j\geq 1}$ be a sequence of centered and
real valued random variables satisfying (\ref{hypoalpha}), (\ref{hypoQ}) and
(\ref{hypogamma}). Let $A$ and $\ell $ be two positive integers such that $%
A2^{-\ell }\geq (1\vee 2c_{0}^{-1})$. Let $M=H^{-1}(\tau (c^{-1/\gamma
_{1}}A))$ and for any $j$, set $\overline{X}_{M}(j)=\varphi _{M}(X_{j})-%
\mathbb{E}\varphi _{M}(X_{j})$. \noindent Then, there exists a subset $%
K_{A}^{(\ell )}$ of $\{1,\dots ,A\}$ with $\mathrm{Card}(K_{A}^{(\ell
)})\geq A/2$, such that for any positive $t\leq \kappa \bigl(A^{\gamma
-1}\wedge (2^{\ell }/A)\bigr)^{\gamma _{1}/\gamma }$, where $\kappa $ is
defined by (\ref{defkappa}), 
\begin{equation}  \label{resultpropinter2}
\log \exp \Bigl(t\sum_{j\in K_{A}^{(\ell )}}\overline{X}_{M}(j)\Bigr)\leq
t^{2}v^{2}A+t^{2}\bigl(\ell (2A)^{1+\frac{\gamma _{1}}{\gamma }}+4A^{\gamma
}(2A)^{\frac{2\gamma _{1}}{\gamma }}\bigr)\exp \Bigl(-\frac{1}{2}\Bigl(\frac{%
c_{1}A}{2^{\ell }}\Bigr)^{\gamma _{1}}\Bigr)\,,
\end{equation}%
with 
\begin{equation}
v^{2}=\sup_{T\geq 1}\sup_{K\subset {\mathbb{N}}^{\ast }}\frac{1}{\mathrm{Card%
}K}\mathrm{Var}\sum_{i\in K}\varphi _{T}(X_{i})  \label{hypovarcont}
\end{equation}%
(the maximum being taken over all nonempty finite sets $K$ of integers).
\end{prop}

\begin{rmk}
\label{compv2V} Notice that $v^2 \leq V$ (the proof is immediate).
\end{rmk}

\noindent \textbf{Proof of Proposition \ref{propinter2}.} The proof is
divided in several steps.

\textit{Step 1. The construction of $K_{A}^{(\ell )}$}. Let $c_{0}$ be
defined by (\ref{defc0}) and $n_{0}=A$. $K_{A}^{(\ell )}$ will be a finite
union of $2^{\ell }$ disjoint sets of consecutive integers with same
cardinal spaced according to a recursive "Cantor"-like construction. We
first define an integer $d_{0}$ as follows: 
\begin{equation*}
d_{0}=\left\{ 
\begin{array}{ll}
\sup \{m\in {2\mathbb{N}}\,,\,m\leq c_{0}n_{0}\} & \text{ if $n_{0}$ is even}
\\ 
\sup \{m\in {2\mathbb{N}}+1\,,\,m\leq c_{0}n_{0}\} & \text{ if $n_{0}$ is odd%
}.%
\end{array}%
\right.
\end{equation*}%
It follows that $n_{0}-d_{0}$ is even. Let $n_{1}=(n_{0}-d_{0})/2$, and
define two sets of integers of cardinal $n_{1}$ separated by a gap of $d_{0}$
integers as follows 
\begin{eqnarray*}
I_{1,1} &=&\{1,\dots ,n_{1}\} \\
I_{1,2} &=&\{n_{1}+d_{0}+1,\dots ,n_{0}\}\,.
\end{eqnarray*}%
We define now the integer $d_{1}$ by 
\begin{equation*}
d_{1}=\left\{ 
\begin{array}{ll}
\sup \{m\in {2\mathbb{N}}\,,\,m\leq c_{0}2^{-(\ell \wedge \frac{1}{\gamma }%
)}n_{0}\} & \text{ if $n_{1}$ is even} \\ 
\sup \{m\in {2\mathbb{N}}+1\,,\,m\leq c_{0}2^{-(\ell \wedge \frac{1}{\gamma }%
)}n_{0}\} & \text{ if $n_{1}$ is odd}.%
\end{array}%
\right.
\end{equation*}%
Noticing that $n_{1}-d_{1}$ is even, we set $n_{2}=(n_{1}-d_{1})/2$, and
define four sets of integers of cardinal $n_{2}$ by 
\begin{eqnarray*}
I_{2,1} &=&\{1,\dots ,n_{2}\} \\
I_{2,2} &=&\{n_{2}+d_{1}+1,\dots ,n_{1}\} \\
I_{2,i+2} &=&(n_{1}+d_{0})+I_{2,i}\text{ for $i=1,2$}\,.
\end{eqnarray*}%
Iterating this procedure $j$ times (for $1\leq j\leq \ell )$, we then get a
finite union of $2^{j}$ sets, $(I_{j,k})_{1\leq k\leq 2^{j}}$, of
consecutive integers, with same cardinal, constructed by induction from $%
(I_{j-1,k})_{1\leq k\leq 2^{j-1}}$ as follows: First, for $1\leq k\leq
2^{j-1}$, we have 
\begin{equation*}
I_{j-1,k}=\{a_{j-1,k},\dots ,b_{j-1,k}\}\,,
\end{equation*}%
where $1+b_{j-1,k}-a_{j-1,k}=n_{j-1}$ and 
\begin{equation*}
1=a_{j-1,1}<b_{j-1,1}<a_{j-1,2}<b_{j-1,2}<\cdots
<a_{j-1,2^{j-1}}<b_{j-1,2^{j-1}}=n_{0}.
\end{equation*}%
Let $n_{j}=2^{-1}(n_{j-1}-d_{j-1})$ and 
\begin{equation*}
d_{j}=\left\{ 
\begin{array}{ll}
\sup \{m\in {2\mathbb{N}}\,,\,m\leq c_{0}2^{-(\ell \wedge \frac{j}{\gamma }%
)}n_{0}\} & \text{ if $n_{j}$ is even} \\ 
\sup \{m\in {2\mathbb{N}}+1\,,\,m\leq c_{0}2^{-(\ell \wedge \frac{j}{\gamma }%
)}n_{0}\} & \text{ if $n_{j}$ is odd}.%
\end{array}%
\right.
\end{equation*}%
Then $I_{j,k}=\{a_{j,k},a_{j,k}+1,\dots ,b_{j,k}\}$, where the double
indexed sequences $(a_{j,k})$ and $(b_{j,k})$ are defined as follows: 
\begin{equation*}
a_{j,2k-1}=a_{j-1,k}\,,\,b_{j,2k}=b_{j-1,k}\,,\,b_{j,2k}-a_{j,2k}+1=n_{j}%
\text{ and }b_{j,2k-1}-a_{j,2k-1}+1=n_{j}\,.
\end{equation*}%
With this selection, we then get that there is exactly $d_{j-1}$ integers
between $I_{j,2k-1}$ and $I_{j,2k}$ for any $1\leq k\leq 2^{j-1}$.

Finally we get 
\begin{equation*}
K_{A}^{(\ell )}=\bigcup_{k=1}^{2^{\ell }}I_{\ell ,k}\,.
\end{equation*}%
Since $\mathrm{Card}(I_{\ell ,k})=n_{\ell }$, for any $1\leq k\leq 2^{\ell }$%
, we get that $\mathrm{Card}(K_{A}^{(\ell )})=2^{\ell }n_{\ell }$. Now
notice that 
\begin{equation*}
A-\mathrm{Card}(K_{A}^{(\ell )})=\sum_{j=0}^{\ell -1}2^{j}d_{j}\leq Ac_{0}%
\Big (\sum_{j\geq 0}2^{j(1-1/\gamma )}+\sum_{j\geq 1}2^{-j}\Big )\leq A/2\,.
\end{equation*}%
Consequently 
\begin{equation*}
A\geq \mathrm{Card}(K_{A}^{(\ell )})\geq A/2\,\text{ and }\,n_{\ell }\leq
A2^{-\ell }\,.
\end{equation*}

The following notation will be useful for the rest of the proof: For any $k$
in $\{0,1,\dots, \ell \}$ and any $j$ in $\{1,\dots, 2^{\ell}\}$ , we set 
\begin{equation}  \label{defKAlj}
K^{(\ell)}_{A, k, j} = \bigcup_{i=(j-1)2^{\ell - k} +1}^{j2^{\ell -
k}}I_{\ell , i} \, .
\end{equation}
Notice that $K_{A}^{( \ell)}= K^{(\ell)}_{A, 0, 1}$ and that for any $k$ in $%
\{0,1,\dots, \ell\}$ 
\begin{equation}  \label{defKAl}
K_{A}^{( \ell)} = \bigcup_{j=1}^{2^{k}}K^{(\ell)}_{A,k, j} \, ,
\end{equation}
where the union is disjoint.

\medskip

In what follows we shall also use the following notation: for any integer $j$
in $[0,\ell ]$, we set 
\begin{equation}
M_{j}=H^{-1}\bigl(\tau (c^{-1/\gamma _{1}}A2^{-(\ell \wedge \frac{j}{\gamma }%
)})\bigr)\,.  \label{defMk}
\end{equation}%
Since $H^{-1}(y)=\big (\log (e/y)\big)^{1/\gamma _{2}}$ for any $y\leq e$,
we get that for any $x\geq 1$, 
\begin{equation}
H^{-1}(\tau (c^{-1/\gamma _{1}}x))\leq \big (1+x^{\gamma _{1}}\big)%
^{1/\gamma _{2}}\leq (2x)^{\gamma _{1}/\gamma _{2}}\,.  \label{majQalpha}
\end{equation}%
Consequently since for any $j$ in $[0,\ell ]$, $A2^{-(\ell \wedge \frac{j}{%
\gamma })}\geq 1$, the following bound is valid: 
\begin{equation}
M_{j}\leq \big (2A2^{-(\ell \wedge \frac{j}{\gamma })}\big )^{\gamma
_{1}/\gamma _{2}}\,.  \label{majMl}
\end{equation}%
For any set of integers $K$ and any positive $M$ we also define 
\begin{equation}
\overline{S}_{M}(K)=\sum_{i\in K}\overline{X}_{M}(i)\,.  \label{defSMK}
\end{equation}

\medskip

\textit{Step 2. Proof of Inequality (\ref{resultpropinter2}) with $%
K_A^{(\ell)}$ defined in step 1}.

Consider the decomposition (\ref{defKAl}), and notice that for any $i=1,2$, $%
\mathrm{Card}(K_{A,1,i}^{(\ell )})\leq A/2$ and 
\begin{equation*}
\tau \bigl(\sigma (X_{i}\,:\,i\in K_{A,1,1}^{(\ell )}),\bar{S}%
_{M_{0}}(K_{A,1,2}^{(\ell )})\bigr)\leq A\tau (d_{0})/2\,.
\end{equation*}%
Since $\overline{X}_{M_{0}}(j)\leq 2M_{0}$, we get that $|\bar{S}%
_{M_{0}}(K_{A,1,i}^{(\ell )})|\leq AM_{0}$. Consequently, by using Lemma 2
from Appendix, we derive that for any positive $t$, 
\begin{equation*}
|\mathbb{E}\exp \big (t\bar{S}_{M_{0}}(K_{A}^{(\ell )})\big )-\prod_{i=1}^{2}%
\mathbb{E}\exp \big (t\bar{S}_{M_{0}}(K_{A,1,i}^{(\ell )})\big )|\leq \frac{%
At}{2}\tau (d_{0})\exp (2tAM_{0})\,.
\end{equation*}%
Since the random variables $\bar{S}_{M_{0}}(K_{A}^{(\ell )})$ and $\bar{S}%
_{M_{0}}(K_{A,1,i}^{(\ell )})$ are centered, their Laplace transform are
greater than one. Hence applying the elementary inequality 
\begin{equation}
|\log x-\log y|\leq |x-y|\text{ for }x\geq 1\text{ and }y\geq 1,
\label{inelog}
\end{equation}%
we get that, for any positive $t$, 
\begin{equation}
|\log \mathbb{E}\exp \big (t\bar{S}_{M_{0}}(K_{A}^{(\ell )})\big )%
-\sum_{i=1}^{2}\log \mathbb{E}\exp \big (t\bar{S}_{M_{0}}(K_{A,1,i}^{(\ell
)})\big )|\leq \frac{At}{2}\tau (d_{0})\exp (2tAM_{0})\,.  \notag  \label{1}
\end{equation}

The next step is to compare $\mathbb{E}\exp \big (t\bar{S}%
_{M_{0}}(K_{A,1,i}^{(\ell )})\big )$ with $\mathbb{E}\exp \big (t\bar{S}%
_{M_{1}}(K_{A,1,i}^{(\ell )})\big )$ for $i=1,2$. The random variables $\bar{%
S}_{M_{0}}(K_{A,1,i}^{(\ell )})$ and $\bar{S}_{M_{1}}(K_{A,1,i}^{(\ell )})$
have values in $[-AM_{0},AM_{0}]$, hence applying the inequality 
\begin{equation}
|e^{tx}-e^{ty}|\leq |t||x-y|(e^{|tx|}\vee e^{|ty|})\,,  \label{AFexp}
\end{equation}%
we obtain that, for any positive $t$, 
\begin{equation*}
\big |\mathbb{E}\exp \big (t\bar{S}_{M_{0}}(K_{A,1,i}^{(\ell )})\big )-%
\mathbb{E}\exp \big (t\bar{S}_{M_{1}}(K_{A,1,i}^{(\ell )})\big )\big |\leq
te^{tAM_{0}}\mathbb{E}\big |\bar{S}_{M_{0}}(K_{A,1,i}^{(\ell )})-\bar{S}%
_{M_{1}}(K_{A,1,i}^{(\ell )})\big |\,.
\end{equation*}%
Notice that 
\begin{equation*}
\mathbb{E}\big |\bar{S}_{M_{0}}(K_{A,1,i}^{(\ell )})-\bar{S}%
_{M_{1}}(K_{A,1,i}^{(\ell )})\big |\leq 2\sum_{j\in K_{A,1,i}^{(\ell )}}%
\mathbb{E}|(\varphi _{M_{0}}-\varphi _{M_{1}})(X_{j})|\,.
\end{equation*}%
Since for all $x\in {\mathbb{R}}$, $|(\varphi _{M_{0}}-\varphi
_{M_{1}})(x)|\leq M_{0}{\ 1\hspace{-1mm}{\mathrm{I}}}_{|x|>M_{1}}$, we get
that 
\begin{equation*}
\mathbb{E}|(\varphi _{M_{0}}-\varphi _{M_{1}})(X_{j})|\leq M_{0}\mathbb{P}%
(|X_{j}|>M_{1})\leq M_{0}\tau (c^{-\frac{1}{\gamma _{1}}}A2^{-(\ell \wedge 
\frac{1}{\gamma })})\,.
\end{equation*}%
Consequently, since $\mathrm{Card}(K_{A,1,i}^{(\ell )})\leq A/2$, for any $%
i=1,2$ and any positive $t$, 
\begin{equation*}
\big |\mathbb{E}\exp \big (t\bar{S}_{M_{0}}(K_{A,1,i}^{(\ell )})\big )-%
\mathbb{E}\exp \big (t\bar{S}_{M_{1}}(K_{A,1,i}^{(\ell )})\big )\big |\leq
tAM_{0}e^{tAM_{0}}\tau (c^{-\frac{1}{\gamma _{1}}}A2^{-(\ell \wedge \frac{1}{%
\gamma })})\,.
\end{equation*}%
Using again the fact that the variables are centered and taking into account
the inequality (\ref{inelog}), we derive that for any $i=1,2$ and any
positive $t$, 
\begin{equation}
\big |\log \mathbb{E}\exp \big (t\bar{S}_{M_{0}}(K_{A,1,i}^{(\ell )})\big )%
-\log \mathbb{E}\exp \big (t\bar{S}_{M_{1}}(K_{A,1,i}^{(\ell )})\big )\big |%
\leq e^{2tAM_{0}}\tau (c^{-\frac{1}{\gamma _{1}}}A2^{-(\ell \wedge \frac{1}{%
\gamma })})\,.  \label{2}
\end{equation}

Now for any $k=1,\dots ,\ell $ and any $i=1,\dots ,2^{k}$, $\mathrm{Card}%
(K_{A,k,i}^{(\ell )})\leq 2^{-k}A$. By iterating the above procedure, we
then get for any $k=1,\dots ,\ell $, and any positive $t$, 
\begin{eqnarray*}
|\sum_{i=1}^{2^{k-1}}\log \mathbb{E}\exp \big (t\bar{S}%
_{M_{k-1}}(K_{A,k-1,i}^{(\ell )})\big ) &-&\sum_{i=1}^{2^{k}}\log \mathbb{E}%
\exp \big (t\bar{S}_{M_{k-1}}(K_{A,k,i}^{(\ell )})\big )| \\
&\leq &2^{k-1}\frac{tA}{2^{k}}\tau (d_{k-1})\exp \big (\frac{2tAM_{k-1}}{%
2^{k-1}}\big )\,,
\end{eqnarray*}%
and for any $i=1,\dots ,2^{k}$, 
\begin{equation*}
|\log \mathbb{E}\exp \big (t\bar{S}_{M_{k-1}}(K_{A,k,i}^{(\ell )})\big )%
-\log \mathbb{E}\exp \big (t\bar{S}_{M_{k}}(K_{A,k,i}^{(\ell )})\big )|\leq
\tau (c^{-\frac{1}{\gamma _{1}}}A2^{-(\ell \wedge \frac{k}{\gamma })})\exp %
\big (\frac{2tAM_{k-1}}{2^{k-1}}\big )\,.
\end{equation*}%
Hence finally, we get that for any $j=1,\dots ,\ell $, and any positive $t$, 
\begin{eqnarray*}
&&|\sum_{i=1}^{2^{j-1}}\log \mathbb{E}\exp \big (t\bar{S}%
_{M_{j-1}}(K_{A,j-1,i}^{(\ell )})\big )-\sum_{i=1}^{2^{j}}\log \mathbb{E}%
\exp \big (t\bar{S}_{M_{j}}(K_{A,j,i}^{(\ell )})\big )| \\
&&\quad \quad \leq \frac{tA}{2}\tau (d_{j-1})\exp
(2tAM_{j-1}2^{1-j})+2^{j}\tau (c^{-\frac{1}{\gamma _{1}}}A2^{-(\ell \wedge 
\frac{j}{\gamma })})\exp (2tAM_{j-1}2^{1-j})\,.
\end{eqnarray*}%
Set 
\begin{equation*}
k_{\ell }=\sup \{j\in \mathbb{N}\,,\,j/\gamma <\ell \}\,,
\end{equation*}%
and notice that $0\leq k_{\ell }\leq \ell -1$. Since $K_{A}^{(\ell
)}=K_{A,0,1}^{(\ell )}$, we then derive that for any positive $t$, 
\begin{eqnarray}
&&|\log \mathbb{E}\exp \big (t\bar{S}_{M_{0}}(K_{A}^{(\ell )})\big )%
-\sum_{i=1}^{2^{k_{\ell }+1}}\log \mathbb{E}\exp \big (t\bar{S}_{M_{k_{\ell
}+1}}(K_{A,k_{\ell }+1,i}^{(\ell )})\big )|  \notag  \label{1step} \\
&&\quad \quad \leq \frac{tA}{2}\sum_{j=0}^{k_{\ell }}\tau (d_{j})\exp \big (%
\frac{2tAM_{j}}{2^{j}}\big )+2\sum_{j=0}^{k_{\ell }-1}2^{j}\tau
(2^{-1/\gamma }c^{-1/\gamma _{1}}A2^{-j/\gamma })\exp \big (\frac{2tAM_{j}}{%
2^{j}}\big )  \notag \\
&&\quad \quad \quad +2^{k_{\ell }+1}\tau (c^{-1/\gamma _{1}}A2^{-\ell })\exp
(2tAM_{k_{\ell }}2^{-k_{\ell }})\,.
\end{eqnarray}
Notice now that for any $i=1,\dots ,2^{k_{\ell }+1}$, $S_{M_{k_{\ell
}+1}}(K_{A,k_{\ell }+1,i}^{(\ell )})$ is a sum of $2^{\ell -k_{\ell }-1}$
blocks, each of size $n_{\ell }$ and bounded by $2M_{k_{\ell }+1}n_{\ell }$.
In addition the blocks are equidistant and there is a gap of size $%
d_{k_{\ell }+1}$ between two blocks. Consequently, by using Lemma 2 along
with Inequality (\ref{inelog}) and the fact that the variables are centered,
we get that 
\begin{eqnarray}
&&|\log \mathbb{E}\exp \big (t\bar{S}_{M_{k_{\ell }+1}}(K_{A,k_{\ell
}+1,i}^{(\ell )})\big )-\sum_{j=(i-1)2^{\ell -k_{\ell }-1}+1}^{i2^{\ell
-k_{\ell }-1}}\log \mathbb{E}\exp \big (t\bar{S}_{M_{k_{\ell }+1}}(I_{\ell
,j})\big )|  \notag  \label{2step} \\
&&\quad \quad \leq tn_{\ell }2^{\ell }2^{-k_{\ell }-1}\tau (d_{k_{\ell
}+1})\exp (2tM_{k_{\ell }+1}n_{\ell }2^{\ell -k_{\ell }-1})\,.
\end{eqnarray}%
Starting from (\ref{1step}) and using (\ref{2step}) together with the fact
that $n_{\ell }\leq A2^{-\ell }$, we obtain: 
\begin{eqnarray}
&&|\log \mathbb{E}\exp \big (t\bar{S}_{M_{0}}(K_{A}^{(\ell )})\big )%
-\sum_{j=1}^{2^{\ell }}\log \mathbb{E}\exp \big (t\bar{S}_{M_{k_{\ell
}+1}}(I_{\ell ,j})\big )|  \notag  \label{3step} \\
&&\quad \quad \leq \frac{tA}{2}\sum_{j=0}^{k_{\ell }}\tau (d_{j})\exp \big (%
\frac{2tAM_{j}}{2^{j}}\big )+2\sum_{j=0}^{k_{\ell }-1}2^{j}\tau
(2^{-1/\gamma }c^{-1/\gamma _{1}}A2^{-j/\gamma })\exp \big (\frac{2tAM_{j}}{%
2^{j}}\big )  \notag \\
&&\quad \quad +2^{k_{\ell }+1}\tau (c^{-1/\gamma _{1}}A2^{-\ell })\exp \big (%
\frac{2tAM_{k_{\ell }}}{2^{k_{\ell }}}\big )+tA\tau (d_{k_{\ell }+1})\exp
(tM_{k_{\ell }+1}A2^{-k_{\ell }})\,.
\end{eqnarray}%
Notice that for any $j=0,\dots ,\ell -1$, we have $d_{j}+1\geq \lbrack
c_{0}A2^{-(\ell \wedge \frac{j}{\gamma })}]$ and $c_{0}A2^{-(\ell \wedge 
\frac{j}{\gamma })}\geq 2$. Whence 
\begin{equation*}
d_{j}\geq (d_{j}+1)/2\geq c_{0}A2^{-(\ell \wedge \frac{j}{\gamma })-2}\,.
\end{equation*}%
Consequently setting $c_{1}=\min (\frac{1}{4}c^{1/\gamma
_{1}}c_{0},2^{-1/\gamma })$ and using (\ref{hypoalpha}), we derive that for
any positive $t$,%
\begin{eqnarray*}
&&|\log \mathbb{E}\exp \big (t\bar{S}_{M_{0}}(K_{A}^{(\ell )})\big )%
-\sum_{j=1}^{2^{\ell }}\log \mathbb{E}\exp \big (t\bar{S}_{M_{k_{\ell
}+1}}(I_{\ell ,j})\big )| \\
&\leq &\frac{tA}{2}\sum_{j=0}^{k_{\ell }}\exp \Big(-\big (c_{1}A2^{-j/\gamma
}\big )^{\gamma _{1}}+\frac{2tAM_{j}}{2^{j}}\Big )+2\sum_{j=0}^{k_{\ell
}-1}2^{j}\exp \Big(-\big (c_{1}A2^{-j/\gamma }\big )^{\gamma _{1}}+\frac{%
2tAM_{j}}{2^{j}}\Big ) \\
&&\quad \quad +2^{k_{\ell }+1}\exp \Big(-\big (A2^{-\ell }\big )^{\gamma
_{1}}+\frac{2tAM_{k_{\ell }}}{2^{k_{\ell }}}\Big )+tA\exp \Big(-\big (%
c_{1}A2^{-\ell }\big )^{\gamma _{1}}+tM_{k_{\ell }+1}A2^{-k_{\ell }}\Big )\,.
\end{eqnarray*}%
By (\ref{majMl}), we get that for any $0\leq j\leq k_{\ell }$, 
\begin{equation*}
2AM_{j}2^{-j}\leq 2^{\gamma _{1}/\gamma }(2^{-j}A)^{\gamma _{1}/\gamma }\,.
\end{equation*}%
In addition, since $k_{\ell }+1\geq \gamma \ell $ and $\gamma <1$, we get
that 
\begin{equation*}
M_{k_{\ell }+1}\leq (2A2^{-\ell })^{\gamma _{1}/\gamma _{2}}\leq
(2A2^{-\gamma \ell })^{\gamma _{1}/\gamma _{2}}\,.
\end{equation*}%
Whence, 
\begin{equation*}
M_{k_{\ell }+1}A2^{-k_{\ell }}=2M_{k_{\ell }+1}A2^{-(k_{\ell }+1)}\leq
2^{\gamma _{1}/\gamma }A^{\gamma _{1}/\gamma }2^{-\gamma _{1}\ell }\,.
\end{equation*}%
In addition, 
\begin{equation*}
2AM_{k_{\ell }}2^{-k_{\ell }}\leq 2^{2\gamma _{1}/\gamma }(A2^{-k_{\ell
}-1})^{\gamma _{1}/\gamma }\leq 2^{2\gamma _{1}/\gamma }A^{\gamma
_{1}/\gamma }2^{-\gamma _{1}\ell }\,.
\end{equation*}%
Hence, if $t\leq c_{2}A^{\gamma _{1}(\gamma -1)/\gamma }$ where $%
c_{2}=2^{-(1+2\gamma _{1}/\gamma )}c_{1}^{\gamma _{1}}$, we derive that 
\begin{eqnarray*}
&&|\log \mathbb{E}\exp \big (t\bar{S}_{M_{0}}(K_{A}^{(\ell )})\big )%
-\sum_{j=1}^{2^{\ell }}\log \mathbb{E}\exp \big (t\bar{S}_{M_{k_{\ell
}+1}}(I_{\ell ,j})\big )| \\
&\leq &\frac{tA}{2}\sum_{j=0}^{k_{\ell }}\exp \Big(-\frac{1}{2}\big (%
c_{1}A2^{-{j}/\gamma }\big )^{\gamma _{1}}\Big )+2\sum_{j=0}^{k_{\ell
}-1}2^{j}\exp \Big(-\frac{1}{2}\big (c_{1}A2^{-{j}/\gamma }\big )^{\gamma
_{1}}\Big ) \\
&&\quad \quad +(2^{k_{\ell }+1}+tA)\exp \bigl(-(c_{1}A2^{-\ell })^{\gamma
_{1}}/2\bigr)\,.
\end{eqnarray*}%
Since $2^{k_{\ell }}\leq 2^{\ell \gamma }\leq A^{\gamma }$, it follows that
for any $t\leq c_{2}A^{\gamma _{1}(\gamma -1)/\gamma }$, 
\begin{equation}
|\log \mathbb{E}\exp \big (t\bar{S}_{M_{0}}(K_{A}^{(\ell )})\big )%
-\sum_{j=1}^{2^{\ell }}\log \mathbb{E}\exp \big (t\bar{S}_{M_{k_{\ell
}+1}}(I_{\ell ,j})\big )|\leq (2\ell tA+4A^{\gamma })\exp \Bigl(-\frac{1}{2}%
\Bigl(\frac{c_{1}A}{2^{\ell }}\Bigr)^{\gamma _{1}}\Bigr)\,.  \label{4step}
\end{equation}

We bound up now the log Laplace transform of each $\bar S_{M_{k_{\ell}+1}}
(I_{\ell, j})$ using the following fact: from l'Hospital rule for
monotonicity (see Pinelis (2002)), the function $x\mapsto g (x)=x^{-2}(e^x -
x - 1)$ is increasing on ${\mathbb{R}}$. Hence, for any centered random
variable $U$ such that $\|U\|_{\infty} \leq M$, and any positive $t$, 
\begin{equation}  \label{psi}
\mathbb{E} \exp ( t U) \leq 1 + t^2 g(tM) \mathbb{E} (U^2) \, .
\end{equation}
Notice that 
\begin{equation*}
\| \bar S_{M_{k_{\ell}+1}} (I_{\ell, j}) \|_{\infty} \leq 2 M_{k_{\ell}+1}
n_{\ell} \leq 2^{\gamma_1/\gamma } (A2^{-\ell} )^{\gamma_1/\gamma} .
\end{equation*}
Since $t \leq 2^{-\gamma_1/\gamma } (2^{\ell}/ A )^{\gamma_1/\gamma} $, by
using (\ref{hypovarcont}), we then get that 
\begin{equation*}
\log \mathbb{E} \exp \big (t \bar S_{M_{k_{\ell}+1}} (I_{\ell, j}) \big ) %
\leq t^2 v ^2 n_{\ell} \, .
\end{equation*}
Consequently, for any $t \leq \kappa \big ( A^{ \gamma_1(\gamma -1)/\gamma}
\wedge (2^{ \ell} / A)^{\gamma_1/\gamma})$, the following inequality holds: 
\begin{equation}  \label{5step}
\log \mathbb{E} \exp \big (t \bar S_{M_{0}} (K_{A}^{( \ell)}) \big ) \leq
t^2 v^2 A + (2 \ell t A + 4A^{\gamma} ) \exp \bigl( -\big ( c_1
A2^{-\ell})^{\gamma_1} / 2 \bigr) \, .
\end{equation}
Notice now that $\| \bar S_{M_0} (K_A^{(\ell)}) \|_{\infty} \leq 2 M_0 A
\leq 2^{\gamma_1/\gamma }A^{\gamma_1/\gamma} $. Hence if $t \leq
2^{-\gamma_1/\gamma }A^{-\gamma_1/\gamma}$, by using (\ref{psi}) together
with (\ref{hypovarcont}), we derive that 
\begin{equation}  \label{6step}
\log \mathbb{E} \exp \big (t \bar S_{ M_0 } (K_A^{( \ell)}) \big ) \leq t^2
v^2 A \, ,
\end{equation}
which proves (\ref{resultpropinter2}) in this case.

Now if $2^{-\gamma_1/\gamma }A^{-\gamma_1/\gamma} \leq t \leq \kappa \bigl( %
A^{ \gamma_1(\gamma -1)/\gamma} \wedge ( 2^{ \ell}/A )^{\gamma_1/\gamma} %
\bigr)$, by using (\ref{5step}), we derive that (\ref{resultpropinter2})
holds, which completes the proof of Proposition 1. $\diamond$

\bigskip

We now bound up the Laplace transform of the sum of truncated random
variables on $[1,A]$. Let 
\begin{equation}
\mu =\bigl(2(2\vee 4c_{0}^{-1})/(1-\gamma )\bigr)^{\frac{2}{1-\gamma }}\text{
and }c_{4}=2^{\gamma _{1}/\gamma }3^{\gamma _{1}/\gamma _{2}}c_{0}^{-\gamma
_{1}/\gamma _{2}}\,,  \label{borneA}
\end{equation}%
where $c_{0}$ is defined in (\ref{defc0}). Define also 
\begin{equation}
\nu =\bigl(c_{4}\bigl(3-2^{(\gamma -1)\frac{\gamma _{1}}{\gamma }}\bigr)%
+\kappa ^{-1}\bigr)^{-1}\bigl(1-2^{(\gamma -1)\frac{\gamma _{1}}{\gamma }}%
\bigr)\,,  \label{defK}
\end{equation}%
where $\kappa $ is defined by (\ref{defkappa}).

\begin{prop}
\label{propinter3} Let $(X_{j})_{j\geq 1}$ be a sequence of centered real
valued random variables satisfying (\ref{hypoalpha}), (\ref{hypoQ}) and (\ref%
{hypogamma}). Let $A$ be an integer. Let $M=H^{-1}(\tau (c^{-1/\gamma
_{1}}A))$ and for any $j$, set $\overline{X}_{M}(j)=\varphi _{M}(X_{j})-%
\mathbb{E}\varphi _{M}(X_{j})$. \noindent Then, if $A\geq \mu $ with $\mu $
defined by (\ref{borneA}), for any positive $t<\nu A^{\gamma _{1}(\gamma
-1)/\gamma }$, where $\nu $ is defined by (\ref{defK}), we get that 
\begin{equation}
\log \mathbb{E}\Bigl(\exp ({\textstyle t\sum_{k=1}^{A}\overline{X}_{M}(k)})%
\Bigr)\leq \frac{AV(A)t^{2}}{1-t\nu ^{-1}A^{\gamma _{1}(1-\gamma )/\gamma }}%
\,,  \label{resultpropinter3}
\end{equation}%
where $V(A)=50v^{2}+\nu _{1}\exp (-\nu _{2}A^{\gamma _{1}(1-\gamma )}(\log
A)^{-\gamma })$ and $\nu _{1}$, $\nu _{2}$ are positive constants depending
only on $c$, $\gamma $ and $\gamma _{1}$, and $v^{2}$ is defined by (\ref%
{hypovarcont}).
\end{prop}

\noindent \textbf{Proof of Proposition \ref{propinter3}.} Let $A_{0}=A$ and $%
X^{(0)}(k)=X_{k}$ for any $k=1,\dots ,A_{0}$. Let $\ell $ be a fixed
positive integer, to be chosen later, which satisfies 
\begin{equation}
A_{0}2^{-\ell }\geq (2\vee 4c_{0}^{-1})\,.  \label{rest1l}
\end{equation}%
Let $K_{A_{0}}^{(\ell )}$ be the discrete Cantor type set as defined from $%
\{1,\dots ,A\}$ in Step 1 of the proof of Proposition \ref{propinter2}. Let $%
A_{1}=A_{0}-\mathrm{Card}K_{A_{0}}^{(\ell )}$ and define for any $k=1,\dots
,A_{1}$, 
\begin{equation*}
X^{(1)}(k)=X_{i_{k}}\text{ where }\{i_{1},\dots ,i_{A_{1}}\}=\{1,\dots
,A\}\setminus K_{A}\,.
\end{equation*}%
Now for $i\geq 1$, let $K_{A_{i}}^{(\ell _{i})}$ be defined from $\{1,\dots
,A_{i}\}$ exactly as $K_{A}^{(\ell )}$ is defined from $\{1,\dots ,A\}$.
Here we impose the following selection of $\ell _{i}$: 
\begin{equation}
\ell _{i}=\inf \{j\in \mathbb{N}\,,\,A_{i}2^{-j}\leq A_{0}2^{-\ell }\,\}\,.
\label{restli}
\end{equation}%
Set $A_{i+1}=A_{i}-\mathrm{Card}K_{A_{i}}^{(\ell _{i})}$ and $\{j_{1},\dots
,j_{A_{i+1}}\}=\{1,\dots ,A_{i+1}\}\setminus K_{A_{i+1}}^{(\ell _{i+1})}$.
Define now 
\begin{equation*}
X^{(i+1)}(k)=X^{(i)}(j_{k})\text{ for }k=1,\dots ,A_{i+1}\,.
\end{equation*}%
%
%
%
Let 
\begin{equation}
m(A)=\inf \{m\in \mathbb{N}\,,\,A_{m}\leq A2^{-\ell }\}\,.  \label{defm(A)}
\end{equation}%
Note that $m(A)\geq 1$, since $A_{0}>A2^{-\ell }$ ($\ell \geq 1$). In
addition, $m(A)\leq \ell $ since for all $i\geq 1$, $A_{i}\leq A2^{-i}$.

Obviously, for any $i=0,\dots ,m(A)-1$, the sequences $(X^{(i+1)}(k))$
satisfy (\ref{hypoalpha}), (\ref{hypoQ}) and (\ref{hypovarcont}) with the
same constants. Now we set $T_{0}=M=H^{-1}(\tau (c^{-1/\gamma _{1}}A_{0}))$,
and for any integer $j=0,\dots ,m(A)$,%
\begin{equation*}
T_{j}=H^{-1}(\tau (c^{-1/\gamma _{1}}A_{j}))\,.
\end{equation*}%
With this definition, we then define for all integers $i$ and $j$, 
\begin{equation*}
X_{T_{j}}^{(i)}(k)=\varphi _{T_{j}}\big (X^{(i)}(k)\big )-\mathbb{E}\varphi
_{T_{j}}\big (X^{(i)}(k)\big )\,.
\end{equation*}%
Notice that by (\ref{hypoalpha}) and (\ref{hypoQ}), we have that for any
integer $j\geq 0$, 
\begin{equation}
T_{j}\leq (2A_{j})^{\gamma _{1}/\gamma _{2}}\,.  \label{boundTl}
\end{equation}%
For any $j=1,\dots ,m(A)$ and $i<j$, define 
\begin{equation*}
Y_{i}=\sum_{k\in K_{A_{i}}^{(\ell _{i})}}X_{T_{i}}^{(i)}(k)\,,\
Z_{i}=\sum_{k=1}^{A_{i}}(X_{T_{i-1}}^{(i)}(k)-X_{T_{i}}^{(i)}(k))\text{ for $%
i>0$, and }R_{j}=\sum_{k=1}^{A_{j}}X_{T_{j-1}}^{(j)}(k)\,.
\end{equation*}%
The following decomposition holds: 
\begin{equation}
\sum_{k=1}^{A_{0}}X_{T_{0}}^{(0)}(k)=\sum_{i=0}^{m(A)-1}Y_{i}+%
\sum_{i=1}^{m(A)-1}Z_{i}+R_{m(A)}\,.  \label{P3prop3}
\end{equation}%
To control the terms in the decomposition (\ref{P3prop3}), we need the
following elementary lemma.

\begin{lma}
\label{compAi} For any $j = 0, \dots, m (A) -1$, $A_{j+1} \geq \frac 13 c_0
A_j $.
\end{lma}

\noindent \textbf{Proof of Lemma \ref{compAi}. } Notice that for any $i$ in $%
[0, m(A)[$, we have $A_{i+1} \geq [c_0 A_i] - 1$. Since $c_0 A_i \geq 2$, we
derive that $[c_0 A_i] - 1 \geq ( [c_0 A_i] + 1)/3 \geq c_0 A_i/3$, which
completes the proof. $\diamond$

\medskip Using (\ref{boundTl}), a useful consequence of Lemma \ref{compAi}
is that for any $j=1, \dots, m(A)$ 
\begin{equation}  \label{secondboundTl}
2A_{j} T_{j-1} \leq c_4 A_j^{\gamma_1 /\gamma}
\end{equation}
where $c_4$ is defined by (\ref{borneA})

\medskip

\noindent \textit{A bound for the Laplace transform of $R_{m (A)}$.}

The random variable $|R_{m(A)}|$ is a.s. bounded by $2A_{m(A)}T_{{m(A)}-1}$.
By using (\ref{secondboundTl}) and (\ref{defm(A)}), we then derive that 
\begin{equation}
\Vert R_{m(A)}\Vert _{\infty }\leq c_{4}(A_{m(A)})^{\gamma _{1}/\gamma }\leq
c_{4}\bigl(A2^{-\ell }\bigr)^{\gamma _{1}/\gamma }\,.  \label{boundRl}
\end{equation}%
Hence, if $t\leq c_{4}^{-1}(2^{\ell }/A)^{\gamma _{1}/\gamma }$, by using (%
\ref{psi}) together with (\ref{hypovarcont}), we obtain 
\begin{equation}
\log \mathbb{E}\bigl(\exp ({\textstyle tR_{m(A)}})\bigr)\leq
t^{2}v^{2}A2^{-\ell }\leq t^{2}(v\sqrt{A})^{2}:=t^{2}\sigma _{1}^{2}\,.
\label{P5prop3}
\end{equation}

\medskip

\noindent \textit{A bound for the Laplace transform of the $Y_i$'s.}

Notice that for any $0 \leq i \leq m(A)-1$, by the definition of $\ell_i$
and (\ref{rest1l}), we get that 
\begin{equation*}
2^{- \ell_i} A_i = 2^{1 - \ell_i} (A_i /2) > 2^{-\ell } (A/2) \geq ( 1 \vee
2c_0^{-1} ) \, .
\end{equation*}
Now, by Proposition 1, we get that for any $i \in [0, m(A)[$ and any $t \leq
\kappa \big (
A_i^{ \gamma -1} \wedge 2^{- \ell_i} A_i \big )^{\gamma_1/\gamma} $ with $%
\kappa$ defined by (\ref{defkappa}), 
\begin{equation*}
\log \mathbb{E} \bigl( e^{t Y_i} \bigr) \leq t^2 \Big ( v \sqrt{A_i} + \big
( \sqrt{\ell_i} (2A_i)^{\frac{1}{2}+ \frac{\gamma_1}{2\gamma}} +
2A_i^{\gamma/2} (2A_i)^{\gamma_1/\gamma} \big ) \exp \big( - \frac 14 \big ( %
c_1 A_i2^{-\ell_i} \big )^{\gamma_1} \big ) \Big )^2 \, .
\end{equation*}
Notice now that $\ell_i \leq \ell \leq A$, $A_i \leq A2^{-i}$ and $2^{-\ell
-1} A \leq 2^{-\ell_i} A_i \leq 2^{-\ell} A$. Taking into account these
bounds and the fact that $\gamma < 1$, we then get that for any $i$ in $[0 ,
m(A)[$ and any $t \leq \kappa (2^i /A)^{ 1- \gamma} \wedge (2^\ell
/A)^{\gamma_1/\gamma} $, 
\begin{equation}  \label{P6prop3}
\log \mathbb{E} \bigl( e^{t Y_i} \bigr) \leq t^2 \Bigl( v \frac{A^{1/2}}{%
2^{i/2}} + \Bigl( 2^{2+ \frac{\gamma_1}{\gamma} } \frac{A^{ 1+ \frac{\gamma_1%
}{\gamma} } }{ (2^i)^{ \frac{\gamma}{2} +\frac{\gamma_1}{2\gamma} } } \Bigr) %
\exp \Bigl( - \frac{c^{\gamma_1}_1}{2^{2+\gamma_1}} \Bigl( \frac{A}{2^\ell} %
\Bigr)^{\gamma_1} \Bigr) \Bigr)^2 := t^2 \sigma_{2,i}^2 \, ,
\end{equation}

\medskip

\noindent \textit{A bound for the Laplace transform of the $Z_i$'s.}

Notice first that for any $1\leq i\leq m(A)-1$, $Z_{i}$ is a centered random
variable, such that 
\begin{equation*}
|Z_{i}|\leq \sum_{k=1}^{A_{i}}\Big (\big |(\varphi _{T_{i-1}}-\varphi
_{T_{i}})X^{(i)}(k)\big |+\mathbb{E}|(\varphi _{T_{i-1}}-\varphi
_{T_{i}})X^{(i)}(k)\big |\Big )\,.
\end{equation*}%
Consequently, using (\ref{secondboundTl}) we get that 
\begin{equation*}
\Vert Z_{i}\Vert _{\infty }\leq 2A_{i}T_{i-1}\leq c_{4}A_{i}^{\gamma
_{1}/\gamma }\,.
\end{equation*}%
In addition, since $|(\varphi _{T_{i-1}}-\varphi _{T_{i}})(x)|\leq
(T_{i-1}-T_{i}){\ 1\hspace{-1mm}{\mathrm{I}}}_{x>T_{i}}$, and the random
random variables $(X^{(i)}(k))$ satisfy (\ref{hypoQ}), by the definition of $%
T_{i}$, we get that 
\begin{equation*}
\mathbb{E}|Z_{i}|^{2}\leq (2A_{i}T_{i-1})^{2}\tau (c^{-1/\gamma
_{1}}A_{i})\leq c_{4}^{2}A_{i}^{2\gamma _{1}/\gamma }\,.
\end{equation*}%
Hence applying (\ref{psi}) to the random variable $Z_{i}$, we get for any
positive $t$, 
\begin{equation*}
\mathbb{E}\exp (tZ_{i})\leq 1+t^{2}g(c_{4}tA_{i}^{\gamma _{1}/\gamma
})c_{4}^{2}A_{i}^{2\gamma _{1}/\gamma }\exp (-A_{i}^{\gamma _{1}})\,.
\end{equation*}%
Hence, since $A_{i}\leq A2^{-i}$, for any positive $t$ satisfying $t\leq
(2c_{4})^{-1}(2^{i}/A)^{\gamma _{1}(1-\gamma )/\gamma }$, we have that 
\begin{equation*}
2tA_{i}T_{i-1}\leq A_{i}^{\gamma _{1}}/2\,.
\end{equation*}%
Since $g(x)\leq e^{x}$ for $x\geq 0$, we infer that for any positive $t$
with $t\leq (2c_{4})^{-1}(2^{i}/A)^{\gamma _{1}(1-\gamma )/\gamma }$, 
\begin{equation*}
\log \mathbb{E}\exp (tZ_{i})\leq c_{4}^{2}t^{2}(2^{-i}A)^{2\gamma
_{1}/\gamma }\exp (-A_{i}^{\gamma _{1}}/2)\,.
\end{equation*}%
By taking into account that for any $1\leq i\leq m(A)-1$, $A_{i}\geq
A_{m(A)-1}>A2^{-\ell }$ (by definition of $m(A)$), it follows that for any $%
i $ in $[1,m(A)[$ and any positive $t$ satisfying $t\leq
(2c_{4})^{-1}(2^{i}/A)^{\gamma _{1}(1-\gamma )/\gamma }$, 
\begin{equation}
\log \mathbb{E}\exp (tZ_{i})\leq t^{2}\bigl(c_{4}(2^{-i}A)^{\gamma
_{1}/\gamma }\exp (-(A2^{-\ell })^{\gamma _{1}}/4)\bigr)^{2}:=t^{2}\sigma
_{3,i}^{2}\,.  \label{P7prop3}
\end{equation}

\noindent \textit{End of the proof.} Let 
\begin{equation*}
C= c_4\Big ( \frac{A}{2^{\ell}} \Big )^{ \gamma_1/\gamma}+ \frac{1}{\kappa}
\sum_{i=0}^{m(A)-1} \Big ( \Big ( \frac{A}{2^i} \Big )^{ 1- \gamma} \vee 
\frac{A}{2^{ \ell}} \Big )^{\gamma_1/\gamma}+ 2c_4 \sum_{i=1}^{m(A)-1} \Big
( \frac{A}{2^i} \Big )^{ \gamma_1(1- \gamma)/\gamma} \, ,
\end{equation*}
and 
\begin{equation*}
\sigma= \sigma_1 + \sum_{i=0}^{m(A)-1} \sigma_{2,i} + \sum_{i=1}^{m(A)-1}
\sigma_{3,i} \, ,
\end{equation*}
where $\sigma_1$, $\sigma_{2,i}$ and $\sigma_{3,i}$ are respectively defined
in (\ref{P5prop3}), (\ref{P6prop3}) and (\ref{P7prop3}).

Notice that $m(A)\leq \ell $, and $\ell \leq 2\log A/\log 2$. We select now $%
\ell $ as follows 
\begin{equation*}
\ell =\inf \{j\in \mathbb{N}\,,\,2^{j}\geq A^{\gamma }(\log A)^{\gamma
/\gamma _{1}}\}\,.
\end{equation*}%
This selection is compatible with (\ref{rest1l}) if 
\begin{equation}
(2\vee 4c_{0}^{-1})(\log A)^{\gamma /\gamma _{1}}\leq A^{1-\gamma }\,.
\label{restrictA}
\end{equation}%
Now we use the fact that for any positive $\delta $ and any positive $u$, $%
\delta \log u\leq u^{\delta }$. Hence if $A\geq 3$, 
\begin{equation*}
(2\vee 4c_{0}^{-1})(\log A)^{\gamma /\gamma _{1}}\leq (2\vee
4c_{0}^{-1})\log A\leq 2(1-\gamma )^{-1}(2\vee 4c_{0}^{-1})A^{(1-\gamma
)/2}\,,
\end{equation*}%
which implies that (\ref{restrictA}) holds as soon as $A\geq \mu $ where $%
\mu $ is defined by (\ref{borneA}). It follows that 
\begin{equation}
C\leq \nu ^{-1}A^{\gamma _{1}(1-\gamma )/\gamma }\,.  \label{boundC}
\end{equation}%
In addition 
\begin{equation*}
\sigma \leq 5v\sqrt{A}+10\times 2^{2\gamma _{1}/\gamma }A^{1+\gamma
_{1}/\gamma }\exp \big(-\frac{c_{1}^{\gamma _{1}}}{2^{2+\gamma _{1}}}%
(A2^{-\ell }\big )^{\gamma _{1}}\big )+c_{4}A^{\gamma _{1}/\gamma }\exp \big
(-\frac{1}{4}(A2^{-\ell })^{\gamma _{1}}\big )\,.
\end{equation*}%
Consequently, since $A2^{-\ell }\geq \frac{1}{2}A^{1-\gamma }(\log
A)^{-\gamma /\gamma _{1}}$, there exists positive constants $\nu _{1}$ and $%
\nu _{2}$ depending only on $c$, $\gamma $ and $\gamma _{1}$ such that 
\begin{equation}
\sigma ^{2}\leq A\bigl(50v^{2}+\nu _{1}\exp (-\nu _{2}A^{\gamma
_{1}(1-\gamma )}(\log A)^{-\gamma })\bigr)=AV(A)\,.  \label{sigmagene}
\end{equation}%
Starting from the decomposition (\ref{P3prop3}) and the bounds (\ref{P5prop3}%
), (\ref{P6prop3}) and (\ref{P7prop3}), we aggregate the contributions of
the terms by using Lemma \ref{breta} given in the appendix. Then, by taking
into account the bounds (\ref{boundC}) and (\ref{sigmagene}), Proposition 2
follows. $\diamond $

\subsection{Proof of Theorem \protect\ref{BTinegacont}}

For any positive $M$ and any positive integer $i$, we set 
\begin{equation*}
\overline X_M (i) = \varphi_M (X_i) - \mathbb{E} \varphi_M (X_i) \, .
\end{equation*}

\noindent $\bullet $ If $\lambda \geq n^{\gamma _{1}/\gamma }$, setting $%
M=\lambda /n$, we have: 
\begin{equation*}
\sum_{i=1}^{n}|\overline{X}_{M}(i)|\leq 2\lambda ,
\end{equation*}%
which ensures that 
\begin{equation*}
\mathbb{P}\Bigl(\sup_{j\leq n}|S_{j}|\geq 3\lambda \Bigr)\leq \mathbb{P}\Big
(\sum_{i=1}^{n}|X_{i}-\overline{X}_{M}(i)|\geq \lambda \Big ).
\end{equation*}%
Now 
\begin{equation*}
\mathbb{P}\Big (\sum_{i=1}^{n}|X_{i}-\overline{X}_{M}(i)|\geq \lambda \Big )%
\leq \frac{1}{\lambda }\sum_{i=1}^{n}\mathbb{E}|X_{i}-\overline{X}%
_{M}(i)|\leq \frac{2n}{\lambda }\int_{M}^{\infty }H(x)dx.
\end{equation*}%
Now recall that $\log H(x)=1-x^{\gamma _{2}}$. It follows that the function $%
x\rightarrow \log (x^{2}H(x))$ is nonincreasing for $x\geq (2/\gamma
_{2})^{1/\gamma _{2}}$. Hence, for $M\geq (2/\gamma _{2})^{1/\gamma _{2}}$, 
\begin{equation*}
\int_{M}^{\infty }H(x)dx\leq M^{2}H(M)\int_{M}^{\infty }\frac{dx}{x^{2}}%
=MH(M).
\end{equation*}%
Whence 
\begin{equation}
\mathbb{P}\Big (\sum_{i=1}^{n}|X_{i}-\overline{X}_{M}(i)|\geq \lambda \Big )%
\leq 2n\lambda ^{-1}MH(M)\ \hbox{ for any }\ M\geq (2/\gamma _{2})^{1/\gamma
_{2}}.  \label{B1decST*}
\end{equation}
Consequently our choice of $M$ together with the fact that $(\lambda
/n)^{\gamma _{2}}\geq \lambda ^{\gamma }$ lead to 
\begin{equation*}
\mathbb{P}\Bigl(\sup_{j\leq n}|S_{j}|\geq 3\lambda \Bigr)\leq 2\exp
(-\lambda ^{\gamma })\,
\end{equation*}%
provided that $\lambda /n\geq (2/\gamma _{2})^{1/\gamma _{2}}$. Now since $%
\lambda \geq n^{\gamma _{1}/\gamma }$, this condition holds if $\lambda \geq
(2/\gamma _{2})^{1/\gamma }$. Consequently for any $\lambda \geq n^{\gamma
_{1}/\gamma }$, we get that 
\begin{equation*}
\mathbb{P}\Bigl(\sup_{j\leq n}|S_{j}|\geq 3\lambda \Bigr)\leq e\exp
(-\lambda ^{\gamma }/C_{1})\,,
\end{equation*}%
as soon as $C_{1}\geq (2/\gamma _{2})^{1/\gamma }$.

\medskip

\noindent $\bullet $ Let $\zeta =\mu \vee (2/\gamma _{2})^{1/\gamma _{1}}$
where $\mu $ is defined by (\ref{borneA}). Assume that $(4\zeta )^{\gamma
_{1}/\gamma }\leq \lambda \leq n^{\gamma _{1}/\gamma }$. Let $p$ be a real
in $[1,\frac{n}{2}]$, to be chosen later on. Let 
\begin{equation*}
A=\Big [\frac{n}{2p}\Big ],\,k=\Big [\frac{n}{2A}\Big ]\text{ and }%
M=H^{-1}(\tau (c^{-\frac{1}{\gamma _{1}}}A))\,.
\end{equation*}%
For any set of natural numbers $K$, denote 
\begin{equation*}
\bar{S}_{M}(K)=\sum_{i\in K}\overline{X}_{M}(i)\,.
\end{equation*}%
For $i$ integer in $[1,2k]$, let $I_{i}=\{1+(i-1)A,\dots ,iA\}$. Let also $%
I_{2k+1}=\{1+2kA,\dots ,n\}$. Set 
\begin{equation*}
\bar{S}_{1}(j)=\sum_{i=1}^{j}\bar{S}_{M}(I_{2i-1})\ \hbox{ and }\bar{S}%
_{2}(j)=\sum_{i=1}^{j}\bar{S}_{M}(I_{2i}).
\end{equation*}%
We then get the following inequality 
\begin{equation}
\sup_{j\leq n}|S_{j}|\leq \sup_{j\leq k+1}|\bar{S}_{1}(j)|+\sup_{j\leq k}|%
\bar{S}_{2}(j)|+2AM+\sum_{i=1}^{n}|X_{i}-\overline{X}_{M}(i)|\,.
\label{decST}
\end{equation}%
Using (\ref{B1decST*}) together with (\ref{hypoalpha}) and our selection of $%
M$, we get for all positive $\lambda $ that 
\begin{equation*}
\mathbb{P}\Big (\sum_{i=1}^{n}|X_{i}-\overline{X}_{M}(i)|\geq \lambda \Big )%
\leq 2n\lambda ^{-1}M\exp (-A^{\gamma _{1}})\ \hbox{ for }\ A\geq (2/\gamma
_{2})^{1/\gamma _{1}}\,.
\end{equation*}%
By using Lemma 5 in Dedecker and Prieur (2004), we get the existence of
independent random variables $(\bar{S}_{M}^{\ast }(I_{2i}))_{1\leq i\leq k}$
with the same distribution as the random variables $\bar{S}_{M}(I_{2i})$
such that 
\begin{equation}
\mathbb{E}|\bar{S}_{M}(I_{2i})-\bar{S}_{M}^{\ast }(I_{2i})|\leq A\tau
(A)\leq A\exp \big (-cA^{\gamma _{1}}\big )\,.  \label{coupY}
\end{equation}%
The same is true for the sequence $(\bar{S}_{M}(I_{2i-1}))_{1\leq i\leq k+1}$%
. Hence for any positive $\lambda $ such that $\lambda \geq 2AM$, 
\begin{eqnarray*}
\mathbb{P}\Bigl(\sup_{j\leq n}|S_{j}|\geq 6\lambda \Bigr) &\leq &\lambda
^{-1}A(2k+1)\exp \big (-cA^{\gamma _{1}}\big )+2n\lambda ^{-1}M\exp
(-A^{\gamma _{1}}) \\
&&+\mathbb{P}\Bigl(\max_{j\leq k+1}\Big|\sum_{i=1}^{j}\bar{S}_{M}^{\ast
}(I_{2i-1})\Big|\geq \lambda \Bigr)+\mathbb{P}\Bigl(\max_{j\leq k}\Big|%
\sum_{i=1}^{j}\bar{S}_{M}^{\ast }(I_{2i})\Big|\geq \lambda \Bigr)\,.
\end{eqnarray*}%
For any positive $t$, due to the independence and since the variables are
centered, $(\exp (t\bar{S}_{M}(I_{2i})))_{i}$ is a submartingale. Hence
Doob's maximal inequality entails that for any positive $t$, 
\begin{equation*}
\mathbb{P}\Bigl(\max_{j\leq k}\sum_{i=1}^{j}\bar{S}_{M}^{\ast }(I_{2i})\geq
\lambda \Bigr)\leq e^{-\lambda t}\prod_{i=1}^{k}\mathbb{E}\Bigl(\exp (t\bar{S%
}_{M}(I_{2i}))\Bigr)\,.
\end{equation*}%
To bound the Laplace transform of each random variable $\bar{S}_{M}(I_{2i})$%
, we apply Proposition \ref{propinter3} to the sequences $(X_{i+s})_{i\in {%
\mathbb{Z}}}$ for suitable values of $s$. Hence we derive that, if $A\geq
\mu $ then for any positive $t$ such that $t<\nu A^{\gamma _{1}(\gamma
-1)/\gamma }$ (where $\nu $ is defined by (\ref{defK})),%
\begin{equation*}
\sum_{i=1}^{k}\log \mathbb{E}\Bigl(\exp (t\bar{S}_{M}(I_{2i}))\Bigr)\leq
Akt^{2}\frac{V(A)}{1-t\nu ^{-1}A^{\gamma _{1}(1-\gamma )/\gamma }}\,.
\end{equation*}%
Obviously the same inequalities hold true for the sums associated to $%
(-X_{i})_{i\in \mathbb{Z}}$. Now some usual computations (see for instance
page 153 in Rio (2000)) lead to 
\begin{equation*}
\mathbb{P}\Bigl(\max_{j\leq k}\Big|\sum_{i=1}^{j}\bar{S}_{M}^{\ast }(I_{2i})%
\Big|\geq \lambda \Bigr)\leq 2\exp \Bigl(-\frac{\lambda ^{2}}{%
4AkV(A)+2\lambda \nu ^{-1}A^{\gamma _{1}(1-\gamma )/\gamma }}\Bigr)\,.
\end{equation*}%
Similarly, we obtain that 
\begin{equation*}
\mathbb{P}\Bigl(\max_{j\leq k+1}\Big|\sum_{i=1}^{j}\bar{S}_{M}^{\ast
}(I_{2i-1})\Big|\geq \lambda \Bigr)\leq 2\exp \Bigl(-\frac{\lambda ^{2}}{%
4A(k+1)V(A)+2\lambda \nu ^{-1}A^{\gamma _{1}(1-\gamma )/\gamma }}\Bigr)\,.
\end{equation*}%
Let now $p=n\lambda ^{-\gamma /\gamma _{1}}$. It follows that $2A\leq
\lambda ^{\gamma /\gamma _{1}}$ and, since $M\leq (2A)^{\gamma _{1}/\gamma
_{2}}$, we obtain $2AM\leq (2A)^{\gamma _{1}/\gamma }\leq \lambda $. Also,
since $\lambda \geq (4\zeta )^{\gamma _{1}/\gamma }$, we have $n\geq 4p$
implying that $A\geq 4^{-1}\lambda ^{\gamma /\gamma _{1}}\geq \zeta $. The
result follows from the previous bounds.

To end the proof, we mention that if $\lambda \leq ( 4 \zeta
)^{\gamma_1/\gamma}$, then 
\begin{equation*}
\mathbb{P} \Bigl( \sup_{ j \leq n} |S_j| \geq \lambda \Bigr) \leq 1 \leq
e\exp \Bigl( - \frac{ \lambda^\gamma}{(4\zeta)^{\gamma_1}} \Bigr ) \, ,
\end{equation*}
which is less than $n\exp ( - \lambda^\gamma /C_1 ) $ as soon as $n \geq 3$
and $C_1 \geq (4\zeta)^{\gamma_1}$. $\diamond$

\subsection{Proof of Remark \protect\ref{remv2}}

\label{prremv2} Setting $W_{i}=\varphi _{M}(X_{i})$ we first bound $\mathrm{%
Cov}(W_{i},W_{i+k}))$. Applying Inequality (4.2) of Proposition 1 in
Dedecker and Doukhan (2003), we derive that, for any positive $k$, 
\begin{equation*}
|\mathrm{Cov(}W_{i},W_{i+k})|\leq 2\int_{0}^{\gamma ({\mathcal{M}}%
_{i},W_{i+k})/2}Q_{|W_{i}|}\circ G_{|W_{i+k}|}(u)du\,
\end{equation*}%
where 
\begin{equation*}
\gamma ({\mathcal{M}}_{i},W_{i+k})=\Vert \mathbb{E}(W_{i+k}|{\mathcal{M}}%
_{i})-\mathbb{E}(W_{i+k})\Vert _{1}\leq \tau (k)\,,
\end{equation*}%
since $x\mapsto \varphi _{M}(x)$ is $1$-Lipschitz. Now for any $j$, $%
Q_{|W_{j}|}\leq Q_{|X_{j}|}\leq Q$, implying that $G_{|W_{j}|}\geq G$, where 
$G$ is the inverse function of $u\rightarrow \int_{0}^{u}Q(v)dv$. Taking $%
j=i $ and $j=i+k$, we get that 
\begin{equation*}
|\mathrm{Cov}(W_{i},W_{i+k})|\leq 2\int_{0}^{\tau (k)/2}Q_{|X_{i}|}\circ
G(u)du\,.
\end{equation*}%
Making the change-of-variables $u=G(v)$ we also have 
\begin{equation}
|\mathrm{Cov}(W_{i},W_{i+k})|\leq 2\int_{0}^{G(\tau
(k)/2)}Q_{X_{i}}(u)Q(u)du\,,  \label{boundcov}
\end{equation}%
proving the remark.

\subsection{Proof of Theorem \protect\ref{thmMDPsubgeo11}}

For any $n\geq 1$, let $T=T_{n}$ where $(T_{n})$ is a sequence of real
numbers greater than $1$ such that $\lim_{n\rightarrow \infty }T_{n}=\infty $%
, that will be specified later. We truncate the variables at the level $%
T_{n} $. So we consider 
\begin{equation*}
X_{i}^{\prime }=\varphi _{T_{n}}(X_{i})-\mathbb{E}\varphi _{T_{n}}(X_{i})%
\text{ and }X_{i}^{\prime \prime }=X_{i}-X_{i}^{\prime }\,.
\end{equation*}%
Let $S_{n}^{\prime }=\sum_{i=1}^{n}X_{i}^{\prime }$ and $S_{n}^{\prime
\prime }=\sum_{i=1}^{n}X_{i}^{\prime \prime }$. To prove the result, by
exponentially equivalence lemma in Dembo and Zeitouni (1998, Theorem 4.2.13.
p130), it suffices to prove that for any $\eta >0$, 
\begin{equation}
\limsup_{n\rightarrow \infty }a_{n}\log \mathbb{P}\big (\frac{\sqrt{a_{n}}}{%
\sigma _{n}}|S_{n}^{\prime \prime }|\geq \eta \big )=-\infty \,,
\label{SnS'n}
\end{equation}%
and 
\begin{equation}
\text{$\{\frac{{1}}{\sigma _{n}}S_{n}^{\prime }\}$ satisfies (\ref{mdpdef})
with the good rate function $I(t)=\frac{t^{2}}{2}$}\,.  \label{mdpST}
\end{equation}

To prove (\ref{SnS'n}), we first notice that $|x-\varphi
_{T}(x)|=(|x|-T)_{+} $. Consequently, if 
\begin{equation*}
W_{i}^{\prime }=X_{i}-\varphi _{T}(X_{i})\,,
\end{equation*}%
then $Q_{|W_{i}^{\prime }|}\leq (Q-T)_{+}$. Hence, denoting by $%
V_{T_{n}}^{\prime \prime }$ the upper bound for the variance of $%
S_{n}^{\prime \prime }$ (corresponding to $V$ for the variance of $S_{n}$)
we have, by Remark \ref{remv2} , 
\begin{equation*}
V_{T_{n}}^{\prime \prime }\leq
\int_{0}^{1}(Q(u)-T_{n})_{+}^{2}du+4\sum_{k>0}\int_{0}^{\tau _{W^{\prime
}}(k)/2}(Q(G_{T_{n}}(v))-T_{n})_{+}dv\,.
\end{equation*}%
where $G_{T}$ is the inverse function of $x\rightarrow
\int_{0}^{x}(Q(u)-T)_{+}du$ and the coefficients $\tau _{W^{\prime }}(k)$
are the $\tau $-mixing coefficients associated to $(W_{i}^{\prime })_{i}$.
Next, since $x\rightarrow x-\varphi _{T}(x)$ is $1$-Lipschitz, we have that $%
\tau _{W^{\prime }}(k)\leq \tau _{X}(k)=\tau (k)$. Moreover, $G_{T}\geq G$,
because $(Q-T)_{+}\leq Q$. Since $Q$ is nonincreasing, it follows that 
\begin{equation*}
V_{T_{n}}^{\prime \prime }\leq
\int_{0}^{1}(Q(u)-T_{n})_{+}^{2}du+4\sum_{k>0}\int_{0}^{\tau
(k)/2}(Q(G(v))-T_{n})_{+}dv.
\end{equation*}%
Hence 
\begin{equation}
\lim_{n\rightarrow +\infty }V_{T_{n}}^{\prime \prime }=0.  \label{majvarT}
\end{equation}%
The sequence $(X_{i}^{\prime \prime })$ satisfies (\ref{hypoalpha}) and we
now prove that it satisfies also (\ref{hypoQ}) for $n$ large enough. With
this aim, we first notice that, since $|\mathbb{E}(X_{i}^{\prime \prime })|=|%
\mathbb{E}(X_{i}^{\prime })|\leq T$, $|X_{i}^{\prime \prime }|\leq |X_{i}|$
if $|X_{i}|\geq T$. Now if $|X_{i}|<T$ then $X_{i}^{\prime \prime }=\mathbb{E%
}(\varphi _{T}(X_{i}))$, and 
\begin{equation*}
\big |\mathbb{E}(\varphi _{T_{n}}(X_{i}))\big |\leq \int_{T_{n}}^{\infty
}H(x)dx<1\text{ for $n$ large enough}.
\end{equation*}%
Then for $t\geq 1$, 
\begin{equation*}
\sup_{i\in \lbrack 1,n]}\mathbb{P}(|X_{i}^{\prime \prime }|\geq t)\leq
H(t)\,,
\end{equation*}%
proving that the sequence $(X_{i}^{\prime \prime })$ satisfies (\ref{hypoQ})
for $n$ large enough. Consequently, for $n$ large enough, we can apply
Theorem \ref{BTinegacont} to the sequence $(X_{i}^{\prime \prime })$, and we
get that, for any $\eta >0$, 
\begin{equation*}
\mathbb{P}\Bigl(\sqrt{\frac{a_{n}}{\sigma _{n}^{2}}}|S_{n}^{\prime \prime
}|\geq \eta \Bigr)\leq n\exp \Bigl(-\frac{\eta ^{\gamma }\sigma _{n}^{\gamma
}}{C_{1}a_{n}^{\frac{\gamma }{2}}}\Bigr)+\exp \Bigl(-\frac{\eta ^{2}\sigma
_{n}^{2}}{C_{2}a_{n}(1+nV_{T_{n}}^{\prime \prime })}\Bigr)+\exp \Bigl(-\frac{%
\eta ^{2}\sigma _{n}^{2}}{C_{3}na_{n}}\exp \Bigl(\frac{\eta ^{\delta }\sigma
_{n}^{\delta }}{C_{4}a_{n}^{\frac{\delta }{2}}}\Bigr)\Bigr)\,,
\end{equation*}%
where $\delta =\gamma (1-\gamma )/2$. This proves (\ref{SnS'n}), since $%
a_{n}\rightarrow 0,$ $a_{n}n^{\gamma /(2-\gamma )}\rightarrow \infty $ , $%
\lim_{n\rightarrow \infty }V_{T_{n}}^{\prime \prime }=0$ and $\lim
\inf_{n\rightarrow \infty }\sigma _{n}^{2}/n>0$.

\medskip

We turn now to the proof of (\ref{mdpST}). Let $p_{n}=[n^{1/(2-\gamma )}]$
and $q_{n}=\delta _{n}p_{n}$ where $\delta _{n}$ is a sequence of integers
tending to zero and such that%
\begin{equation*}
\delta _{n}^{\gamma _{1}}n^{\gamma _{1}/(2-\gamma )}/\log n\rightarrow
\infty \,\text{ and }\,\delta _{n}^{\gamma _{1}}a_{n}n^{\gamma
_{1}/(2-\gamma )}\rightarrow \infty \,
\end{equation*}%
(this is always possible since $\gamma _{1}\geq \gamma $ and by assumption $%
a_{n}n^{\gamma /(2-\gamma )}\rightarrow \infty $). Let now $%
m_{n}=[n/(p_{n}+q_{n})]$. We divide the variables $\{X_{i}^{\prime }\}$ in
big blocks of size $p_{n}$ and small blocks of size $q_{n}$, in the
following way: Let us set for all $1\leq j\leq m_{n}$, 
\begin{equation*}
Y_{j,n}=\sum_{i=(j-1)(p_{n}+q_{n})+1}^{(j-1)(p_{n}+q_{n})+p_{n}}X_{i}^{%
\prime }\quad \mbox{ and }\quad
Z_{j,n}=\sum_{i=(j-1)(p_{n}+q_{n})+p_{n}+1}^{j\,(p_{n}+q_{n})}X_{i}^{\prime
}\,.
\end{equation*}

Then we have the following decomposition: 
\begin{equation}  \label{dec1Sn'}
S^{\prime }_{n} = \sum_{j=1}^{m_{n}}Y_{j,n}+ \sum_{j=1}^{m_{n}} Z_{j,n} +
\sum_{i=m_n (p_n +q_n) + 1}^{n}X_i^{\prime }\, .
\end{equation}

For any $j=1,\dots ,m_{n}$, let now 
\begin{equation*}
I(n,j)=\{(j-1)(p_{n}+q_{n})+1,\dots ,(j-1)(p_{n}+q_{n})+p_{n}\}\,.
\end{equation*}%
These intervals are of cardinal $p_{n}$. Let 
\begin{equation*}
\ell _{n}=\inf \{k\in \mathbb{N}^{\ast },2^{k}\geq \varepsilon
_{n}^{-1}p_{n}^{\gamma /2}a_{n}^{-1/2}\}\,,
\end{equation*}%
where $\varepsilon _{n}$ a sequence of positive numbers tending to zero and
satisfying 
\begin{equation}
\varepsilon _{n}^{2}a_{n}n^{\gamma /(2-\gamma )}\rightarrow \infty .
\label{choiceepsi}
\end{equation}%
For each $j\in \{1,\dots ,m_{n}\}$, we construct discrete Cantor sets, $%
K_{I(n,j)}^{(\ell _{n})}$, as described in the proof of Proposition \ref%
{propinter2} with $A=p_{n}$, $\ell =\ell _{n}$, and the following selection
of $c_{0}$, 
\begin{equation*}
c_{0}=\frac{\varepsilon _{n}}{1+\varepsilon _{n}}\frac{2^{(1-\gamma )/\gamma
}-1}{2^{1/\gamma }-1}\,.
\end{equation*}%
Notice that clearly with the selections of $p_{n}$ and $\ell _{n}$, $%
p_{n}2^{-\ell _{n}}\rightarrow \infty $. In addition with the selection of $%
c_{0}$ we get that for any $1\leq j\leq m_{n}$, 
\begin{equation*}
\mathrm{Card}(K_{I(n,j)}^{(\ell _{n})})^{c}\leq \frac{\varepsilon _{n}p_{n}}{%
1+\varepsilon _{n}}
\end{equation*}%
and 
\begin{equation*}
K_{I(n,j)}^{(\ell _{n})}=\bigcup_{i=1}^{2^{\ell _{n}}}I_{\ell
_{n},i}(p_{n},j)\,,
\end{equation*}%
where the $I_{\ell _{n},i}(p_{n},j)$ are disjoint sets of consecutive
integers, each of same cardinal such that 
\begin{equation}
\frac{p_{n}}{2^{\ell _{n}}(1+\varepsilon _{n})}\leq \mathrm{Card}I_{\ell
_{n},i}(p_{n},j)\leq \frac{p_{n}}{2^{\ell _{n}}}\,.  \label{majcard}
\end{equation}%
With this notation, we derive that 
\begin{equation}
\sum_{j=1}^{m_{n}}Y_{j,n}=\sum_{j=1}^{m_{n}}S^{\prime }\big (%
K_{I(n,j)}^{(\ell _{n})}\big )+\sum_{j=1}^{m_{n}}S^{\prime }\big (%
(K_{I(n,j)}^{(\ell _{n})})^{c}\big )\,.  \label{dec2Sn'}
\end{equation}%
Combining (\ref{dec1Sn'}) with (\ref{dec2Sn'}), we can rewrite $%
S_{n}^{\prime }$ as follows 
\begin{equation}  \label{construction}
S_{n}^{\prime }=\sum_{j=1}^{m_{n}}S^{\prime }\big (K_{I(n,j)}^{(\ell _{n})}%
\big )+\sum_{k=1}^{r_{n}}\tilde{X}_{i}\,,
\end{equation}%
where $r_{n}=n-m_{n}\mathrm{Card}K_{I(n,1)}^{(\ell _{n})}$ and the $\tilde{X}%
_{i}$ are obtained from the $X_{i}^{\prime }$ and satisfied (\ref{hypoalpha}%
) and (\ref{hypovarcont}) with the same constants. Since $r_{n}=o(n)$,
applying Theorem \ref{BTinegacont} and using the fact that $\lim
\inf_{n\rightarrow \infty }\sigma _{n}^{2}/n>0$, we get that for any $\eta
>0 $, 
\begin{equation}
\limsup_{n\rightarrow \infty }a_{n}\log \mathbb{P}\big (\frac{\sqrt{a_{n}}}{%
\sigma_n}\sum_{k=1}^{r_{n}}\tilde{X}_{i}\geq \eta \big )=-\infty \,.
\label{negligiblereste}
\end{equation}%
Hence to prove (\ref{mdpST}), it suffices to prove that 
\begin{equation}
\text{ $\{\sigma^{-1}_{n}\sum_{j=1}^{m_{n}}S^{\prime }\big (%
K_{I(n,j)}^{(\ell _{n})}\big )\}$ satisfies (\ref{mdpdef}) with the good
rate function $I(t)=t^2/2$}.  \label{mdpST*}
\end{equation}%
With this aim, we choose now $T_n =\varepsilon _{n}^{-1/2}$ where $%
\varepsilon _n$ is defined by (\ref{choiceepsi}).

By using Lemma 5 in Dedecker and Prieur (2004), we get the existence of
independent random variables $\big (S^{\ast }(K_{I(n,j)}^{(\ell _{n})})\big )%
_{1\leq j\leq m_{n}}$ with the same distribution as the random variables $%
S^{\prime }(K_{I(n,j)}^{(\ell _{n})})$ such that 
\begin{equation*}
\sum_{j=1}^{m_{n}}{\mathbb{E}}|S^{\prime }(K_{I(n,j)}^{(\ell _{n})})-S^{\ast
}(K_{I(n,j)}^{(\ell _{n})})|\leq \tau (q_{n})\sum_{j=1}^{m_{n}}\mathrm{Card}%
K_{I(n,j)}^{(\ell _{n})}\,.
\end{equation*}%
Consequently, since $\sum_{j=1}^{m_{n}}\mathrm{Card}K_{I(n,j)}^{(\ell
_{n})}\leq n$, we derive that for any $\eta >0$, 
\begin{equation*}
a_{n}\log \mathbb{P}\big (\frac{\sqrt{a_{n}}}{\sigma_n}%
|\sum_{j=1}^{m_{n}}(S^{\prime }(K_{I(n,j)}^{(\ell _{n})})-S^{\ast
}(K_{I(n,j)}^{(\ell _{n})}))|\geq \eta \big )\leq a_{n}\log \Big (\eta
^{-1}\sigma_n^{-1}n \sqrt{a_{n}}\exp \big (-c\delta _{n}^{\gamma
_{1}}n^{\gamma _{1}/(2-\gamma )}\big )\Big )\,,
\end{equation*}%
which tends to $-\infty $ by the fact that $\liminf_n \sigma _{n}^{2}/n>0$
and the selection of $\delta _{n}$. Hence the proof of the MDP for $%
\{\sigma^{-1} _{n}\sum_{j=1}^{m_{n}}S^{\prime }\big (K_{I(n,j)}^{(\ell _{n})}%
\big
)\}$ is reduced to proving the MDP for $\{\sigma_{n}^{-1}%
\sum_{j=1}^{m_{n}}S^{\ast }\big (K_{I(n,j)}^{(\ell _{n})}\big )\}$. By Ellis
Theorem, to prove (\ref{mdpST*}) it remains then to show that, for any real $%
t$, 
\begin{equation}
a_{n}\sum_{j=1}^{m_{n}}\log \mathbb{E}\exp \Big (tS^{\prime }\big (%
K_{I(n,j)}^{(\ell _{n})}\big )/\sqrt{a_n \sigma _{n}^2}\Big )\rightarrow 
\frac{t^{2}}{2}\text{ as }n\rightarrow \infty \,.  \label{ellisindsg}
\end{equation}%
As in the proof of Proposition \ref{propinter2}, we decorrelate step by
step. Using Lemma \ref{lmainegatau} and taking into account the fact that
the variables are centered together with the inequality (\ref{inelog}), we
obtain, proceeding as in the proof of Proposition \ref{propinter2}, that for
any real $t$, 
\begin{eqnarray*}
&&\Bigl |\sum_{j=1}^{m_{n}}\log \mathbb{E}\exp \Big (tS^{\prime }\big (%
K_{I(n,j)}^{(\ell _{n})}\big )/\sqrt{a_n \sigma _{n}^2}\Big )%
-\sum_{j=1}^{m_{n}}\sum_{i=1}^{2^{\ell _{n}}}\log \mathbb{E}\exp \Big (%
tS^{\prime }\big (I_{\ell _{n},i}(p_{n},j)\big )/\sqrt{a_n \sigma _{n}^2}%
\Big )\Big | \\
&&\quad \leq \frac{|t|m_{n}p_{n}}{\sqrt{a_n \sigma _{n}^2}}\Bigl(\exp \Bigl(%
-c\frac{c_{0}^{\gamma _{1}}}{4}\frac{p_{n}^{\gamma _{1}}}{2^{\ell _{n}\gamma
_{1}}}+2\frac{|t|}{\sqrt{\varepsilon _{n}a_n \sigma _{n}^2}}\frac{p_{n}}{%
2^{\gamma \ell _{n}}}\Bigr)+\sum_{j=0}^{k_{\ell _{n}}}\exp \Bigl(-c\frac{%
c_{0}^{\gamma _{1}}}{4}\frac{p_{n}^{\gamma _{1}}}{2^{j\gamma _{1}/\gamma }}+2%
\frac{|t|}{\sqrt{\varepsilon _{n}a_n \sigma _{n}^2}}\frac{p_{n}}{2^{j}}\Bigr)%
\Bigr)\,,
\end{eqnarray*}%
where $k_{{\ell _{n}}}=\sup \{j\in \mathbb{N}\,,\,j/\gamma <\ell _{n}\}$. By
the selection of $p_{n}$ and $\ell _{n}$, and since $\lim \inf_{n\rightarrow
\infty }\sigma _{n}^{2}/n>0$ and $\varepsilon _{n}^{2}a_{n}n^{\gamma
/(2-\gamma )}\rightarrow \infty $, we derive that for $n$ large enough,
there exists positive constants $K_{1}$ and $K_{2}$ depending on $c$, $%
\gamma $ and $\gamma _{1}$ such that 
\begin{eqnarray}
&&a_{n}\Bigl |\sum_{j=1}^{m_{n}}\log \mathbb{E}\exp \Big (tS^{\prime }\big (%
K_{I(n,j)}^{(\ell _{n})}\big )/\sqrt{a_n \sigma _{n}^2}\Big )%
-\sum_{j=1}^{m_{n}}\sum_{i=1}^{2^{\ell _{n}}}\log \mathbb{E}\exp \Big (%
tS^{\prime }\big (I_{\ell _{n},i}(p_{n},j)\big )/\sqrt{a_n \sigma _{n}^2}%
\Big )\Big |  \notag  \label{dec2proofsg} \\
&\leq &K_{1}|t|\sqrt{a_{n}n}\log (n)\exp \Bigl(-K_{2}\big (\varepsilon
_{n}^{2}a_{n}n^{\gamma /(2-\gamma )}\big )^{\gamma /2}n^{\gamma (1-\gamma
)/(2-\gamma )}\Bigr)\,,
\end{eqnarray}%
which converges to zero by the selection of $\varepsilon _{n}$.

Hence (\ref{ellisindsg}) holds if we prove that for any real $t$ 
\begin{equation}
a_{n}\sum_{j=1}^{m_{n}}\sum_{k=1}^{2^{\ell _{n}}}\log \mathbb{E}\exp \Big (%
tS^{\prime }\big (I_{\ell _{n},i}(p_{n},j)/\sqrt{a_{n}\sigma _{n}^{2}}\Big )%
\rightarrow \frac{t^{2}}{2}\text{ as }n\rightarrow \infty \,.
\label{ellisindsg2}
\end{equation}%
With this aim, we first notice that, by the selection of $\ell _{n}$ and the
fact that $\varepsilon _{n}\rightarrow 0$, 
\begin{equation}
\Vert S^{\prime }\big (I_{\ell _{n},i}(p_{n},j)\Vert _{\infty }\leq
2T_{n}2^{-\ell _{n}}p_{n}=o(\sqrt{na_{n}})=o(\sqrt{\sigma _{n}^{2}a_{n}})\,.
\label{maj1arc}
\end{equation}%
In addition, since $\lim_{n}V_{T_{n}}^{\prime \prime }=0$ and the fact that $%
\lim \inf_{n}\sigma _{n}^{2}/n>0$, we have $\lim_{n}\sigma _{n}^{-2}\mathrm{%
Var}S_{n}^{\prime }=1$. Notice that by (\ref{construction}) and the fact
that $r_{n}=o(n)$, 
\begin{equation*}
\mathrm{Var}S_{n}^{\prime }=\mathbb{E}\Bigl(\sum_{j=1}^{m_{n}}%
\sum_{i=1}^{2^{\ell _{n}}}{S^{\prime }}(I_{\ell _{n},i}(p_{n},j))\Bigr)%
^{2}+o(n)\text{ as }n\rightarrow \infty \,.
\end{equation*}%
Also, a straightforward computation as in Remark \ref{remv2} shows that
under (\ref{hypoalpha}) and (\ref{hypoQ}) we have%
\begin{equation*}
\mathbb{E}\Bigl(\sum_{j=1}^{m_{n}}\sum_{i=1}^{2^{\ell _{n}}}\big ({\
S^{\prime }}\big (I_{\ell _{n},i}(p_{n},j))\Bigr)^{2}=\sum_{j=1}^{m_{n}}%
\sum_{i=1}^{2^{\ell _{n}}}\mathbb{E}\big ({S^{\prime }}^{2}\big (I_{\ell
_{n},i}(p_{n},j)\big )\big )+o(n)\text{ as }n\rightarrow \infty \,.
\end{equation*}%
Hence 
\begin{equation}
\lim_{n\rightarrow \infty }\left( \sigma _{n}\right)
^{-2}\sum_{j=1}^{m_{n}}\sum_{i=1}^{2^{\ell _{n}}}\mathbb{E}\big ({S^{\prime }%
}^{2}\big (I_{\ell _{n},i}(p_{n},j)\big )\big )=1\,.  \label{maj2arc}
\end{equation}%
Consequently (\ref{ellisindsg2}) holds by taking into account (\ref{maj1arc}%
) and (\ref{maj2arc}) and by using Lemma 2.3 in Arcones (2003). $\diamond $

\subsection{ Proof of Corollary \protect\ref{thmMDPsubgeo}}

We have to show that 
\begin{equation}
\lim_{n\rightarrow \infty }\mathrm{Var}(S_{n})/n=\sigma ^{2}>0\,.
\label{limvar}
\end{equation}%
Proceeding as in the proof of Remark \ref{remv2}, we get that for any
positive $k$, 
\begin{equation*}
|\mathrm{Cov}(X_{0},X_{k})|\leq 2\int_{0}^{G(\tau (k)/2)}Q^{2}(u)du\,,
\end{equation*}%
which combined with (\ref{hypoalpha}) and (\ref{hypoQ}) imply that $%
\sum_{k>0}k|\mathrm{Cov}(X_{0},X_{k})|<\infty $. This condition together
with the fact that $\mathrm{Var}(S_{n})\rightarrow \infty $ entails (\ref%
{limvar}) (see Lemma 1 in Bradley (1997)).

\section{Appendix}

\setcounter{equation}{0}

We first give the following decoupling inequality.

\begin{lma}
\label{lmainegatau} Let $Y_{1}$, ..., $Y_{p}$ be real-valued random
variables each a.s. bounded by $M$. For every $i\in \lbrack 1,p]$, let ${%
\mathcal{M}}_{i}=\sigma (Y_{1},...,Y_{i})$ and for $i\geq 2$, let $%
Y_{i}^{\ast }$ be a random variable independent of ${\mathcal{M}}_{i-1}$ and
distributed as $Y_{i}$. Then for any real $t$, 
\begin{equation*}
|\mathbb{E}\exp \Big (t\sum_{i=1}^{p}Y_{i}\Big )-\prod_{i=1}^{p}\mathbb{E}%
\exp (tY_{i})|\leq |t|\exp (|t|Mp)\sum_{i=2}^{p}\mathbb{E}|Y_{i}-Y_{i}^{\ast
}|\,.
\end{equation*}%
In particular, we have for any real $t$, 
\begin{equation*}
|\mathbb{E}\exp \Big (t\sum_{i=1}^{p}Y_{i}\Big )-\prod_{i=1}^{p}\mathbb{E}%
\exp (tY_{i})|\leq |t|\exp (|t|Mp)\sum_{i=2}^{p}\tau (\sigma (Y_{1},\dots
,Y_{i-1}),Y_{i})\,,
\end{equation*}%
where $\tau $ is defined by (\ref{deftau1}).
\end{lma}

\noindent \textbf{Proof of Lemma \ref{lmainegatau}}. Set $%
U_{k}=Y_{1}+Y_{2}+\cdots +Y_{k}$. We first notice that 
\begin{equation}  \label{P1lmatau}
\mathbb{E}\bigl(e^{tU_{p}}\bigr)-\prod_{i=1}^{p}\mathbb{E}\bigl(e^{tY_{i}}%
\bigr)=\sum_{k=2}^{p}\Bigl(\mathbb{E}\bigl(e^{tU_{k}}\bigr)-\mathbb{E}\bigl(%
e^{tU_{k-1}}\bigr)\mathbb{E}\bigl(e^{tY_{k}}\bigr)\Bigr)\prod_{i=k+1}^{p}%
\mathbb{E}\bigl(e^{tY_{i}}\bigr)
\end{equation}%
with the convention that the product from $p+1$ to $p$ has value $1$. Now 
\begin{equation*}
|\mathbb{E}\exp (tU_{k})-\mathbb{E}\exp (tU_{k-1})\mathbb{E}\exp
(tY_{k})|\leq \Vert \exp (tU_{k-1})\Vert _{\infty }\Vert \mathbb{E}\bigl(%
e^{tY_{k}}-e^{tY_{k}^{\ast }}|{{\mathcal{M}}_{k-1}}\bigr)\Vert _{1}\,.
\end{equation*}%
Using (\ref{AFexp}) we then derive that 
\begin{equation}  \label{Plmatau}
|\mathbb{E}\exp (tU_{k})-\mathbb{E}\exp (tU_{k-1})\mathbb{E}\exp
(tY_{k})|\leq |t|\exp (|t|kM)\Vert Y_{k}-Y_{k}^{\ast }\Vert _{1}\,.
\end{equation}
Since the variables are bounded by $M$, starting from (\ref{P1lmatau}) and
using (\ref{Plmatau}), the result follows.$\diamond $

\medskip

One of the tools we use repeatedly is the technical lemma below, which
provides bounds for the log-Laplace transform of any sum of real-valued
random variables.

\begin{lma}
\label{breta} Let $Z_0 , Z_1 , \ldots $ be a sequence of real valued random
variables. Assume that there exist positive constants $\sigma_0 , \sigma_1 ,
\ldots$ and $c_0, c_1 , \ldots $ such that, for any positive $i$ and any $t$
in $[0, 1/c_i[$, 
\begin{equation*}
\log \mathbb{E} \exp (tZ_i) \leq (\sigma_i t)^2 / (1-c_it) \, .
\end{equation*}
Then, for any positive $n$ and any $t$ in $[0, 1/(c_0 + c_1 + \cdots + c_n)[$%
, 
\begin{equation*}
\log \mathbb{E} \exp (t (Z_0 + Z_1 + \cdots + Z_n) ) \leq (\sigma t)^2 /
(1-Ct) ,
\end{equation*}
where $\sigma = \sigma_0 + \sigma_1 + \cdots + \sigma_n $ and $C= c_0 + c_1
+ \cdots + c_n$.
\end{lma}

\noindent \textbf{Proof of Lemma \ref{breta}}. Lemma \ref{breta} follows
from the case $n=1$ by induction on $n$. Let $L$ be the log-Laplace of $%
Z_{0}+Z_{1}$. Define the functions $\gamma _{i}$ by 
\begin{equation*}
\gamma _{i}(t)=(\sigma _{i}t)^{2}/(1-c_{i}t)\ \hbox{ for }t\in \lbrack
0,1/c_{i}[\ \hbox{ and }\gamma _{i}(t)=+\infty \ \hbox{ for }t\geq 1/c_{i}.
\end{equation*}%
For $u$ in $]0,1[$, let $\gamma _{u}(t)=u\gamma _{1}(t/u)+(1-u)\gamma
_{0}(t/(1-u))$. From the H\"{o}lder inequality applied with $p=1/u$ and $%
q=1/(1-u)$, we get that $L(t)\leq \gamma _{u}(t)$ for any nonnegative $t$.
Now, for $t$ in $[0,1/C[$, choose $u=(\sigma _{1}/\sigma )(1-Ct)+c_{1}t$
(here $C=c_{0}+c_{1}$ and $\sigma =\sigma _{0}+\sigma _{1}$). With this
choice $1-u=(\sigma _{0}/\sigma )(1-Ct)+c_{0}t$, so that $u$ belongs to $%
]0,1[$ and 
\begin{equation*}
L(t)\leq \gamma _{u}(t)=(\sigma t)^{2}/(1-Ct),
\end{equation*}%
which completes the proof of Lemma \ref{breta}. $\diamond $


\bigskip

\end{document}